\documentclass[10pt]{amsart}
\usepackage{amsmath,amsthm,amscd,amsfonts,amssymb}
\usepackage[all]{xy}

\begin{document}
\newtheorem{defn}{Definition}[section]
\theoremstyle{definition}
\newtheorem{thm}[defn]{Theorem}
\newtheorem{propo}[defn]{Proposition}
\newtheorem{cor}[defn]{Corollary}
\newtheorem{lem}[defn]{Lemma}
\theoremstyle{remark}
\newtheorem{rem}{Remark}

\renewcommand\o{{{\cal O}}}
\newcommand\s{\sigma}
\newcommand\w{\widehat}
\newcommand\cal{\mathcal}
\newcommand\limi{\varinjlim}
\newcommand\limp{\varprojlim}
\newcommand\te{\text}

\title [ Zeta Correspondences in Rank-$n$ ]{   Zeta Correspondences in Rank-$n$}

\author[A. \'Alvarez]{A. \'Alvarez${}^*$}

\address{Departamento de Matem\'aticas  \\
Universidad de Salamanca \\ Plaza de la Merced 1-4.
Salamanca (37008). Spain.}

\thanks{$*$  Departamento de Matem\'aticas.
Universidad de Salamanca. Spain.}

\maketitle

\tableofcontents


\section{Introduction.}

The purpose of this work is to study in the case of function fields and
rank-$n$   certain correspondences associated with
Jacobi sums   \cite{C}. These correspondences are called zeta
correspondences and have been widely studied \cite{An2},
\cite{An3},
\cite{An4}, \cite{AnDP}, \cite{Si}. Their importance is the relationship
with Stickelberger's Theorem and $L$-series
evaluators, \cite{An2}, \cite{An3}, \cite{An4}, \cite{AnDP}. In work
\cite{C} they are used to prove that the Frobenius
map is essentially a Gauss or Jacobi sum for Fermat and Artin-Schreier curves.

In rank $1$, these zeta correspondences are closely related to the theta
divisor defined over the Jacobian, \cite{An2}. The
strategy in works \cite{An2}, \cite{An3} is to study this subject in the
setting of Drinfeld modules. Here we shall
follow the same strategy but considering ${\Bbb F}_q[t]$-Drinfeld modules
(ellliptic sheaves) with level structures over
an effective divisor $D$ on $Spec({\Bbb F}_q[t])$. In this way, the moduli of
Drinfeld modules (ellliptic sheaves) with level structures are given by the
coefficients of a polynomial. These
coefficients are torsion elements for a universal Drinfeld module
(c.f \cite{Al}): In this setting,  zeta correspondences are given
explicitly in terms of these elements.
 To study rank-$n$ zeta correspondences, we have taken into account the
 generalized theta function considered in \cite{An3} and  the morphism
given in 4. \cite{An2}, but considering the
generalized genus
\cite{Se} pg.86. Firstly, we shall introduce   zeta subschemes over the
moduli of vector bundles with
level structures as a generalization of the $n^2$-generalized theta
divisor. We recall that the $n^2$-generalized
theta divisor over a curve $C$ of genus $g$ is formed by the vector bundles
of rank $n^2$ and degree $n^2(g-1)$ without
global sections. The restriction via a certain morphism, similar to that of
4. \cite{An2}, of  these zeta subschemes to
the moduli of rank-$n$ Drinfeld modules (ellliptic sheaves) will be the
zeta correspondences that we wish to study.

In this way we shall explicitly obtain, for   ${\Bbb P}_1$, the zeta
correspondences as the zero
locus of certain functions. In the rank-$n$ case, they are the zero locus
of $n^2$ functions given by a $n\times
n$-matrix. These functions involve the torsion elements of a rank-$n$
universal ${\Bbb F}_q[t]$-Drinfeld module. To
obtain these results we shall prove a similar Lemma    to "$\chi =0
\Rightarrow h^0=h^1=0$" Lemma 3.3.1, \cite{An2}
for level structures  in the rank-$n$ case and
${\Bbb P}_1$. We can prove this Lemma  in our case because for semistable
bundles over  ${\Bbb P}_1$ the generalized theta
divisor is the empty subset.  For rank-$1$, as in \cite{An2}, we shall see a
relationship  of these zeta correspondences with the theta divisor for
level structures. One observes that instead of  the
zero locus of the $n\times n$-matrix considered one could take the zero
locus of its determinant and would thus obtain a
$1$-correspondence that involves the zeta correspondence   considered before.
The examples that we provide in the rank-$2$ case give us zeta
correspondences whose irreducible components are
essentially the same as in the rank-$1$ case. They are given by graphs of
automorphisms given by central matrices and the
Frobenius morphism.

Above we have
considered Drinfeld modules with level structures over $\infty +D$,
$\infty$ being the pole of the Drinfeld modules.
However if we consider another rational point $p$ different from $\infty$,
we can extend the above results to an
arbitrary proper and smooth curve $C$ over a finite field. We  shall obtain
zeta subschemes for rank-$n$ shtukas; these
zeta subschemes would be formed  by pairs of isogenous shtukas although to
complete the result it would be necessary to
prove  the "$\chi =0 \Rightarrow h^0=h^1=0$"-Lemma for rank $n$. This
result tells us, via the immersion of $A$-Drinfeld modules in  rank-$n$
shtukas, that for ${\Bbb P}_1$  zeta
correspondences ($p$ instead of
$\infty$)  are given by pairs of isogenous $A$-Drinfeld modules.

In the case of a general curve,   these zeta subschemes are again given,
now locally, by the zero locus of
$n^2$-functions. We shall check that these functions are the restriction of
hyperplane  sections via a Plucker
morphism. In this way, these zeta subschemes are the intersection of
$n^2$-generalized theta divisors.

I believe that it would be very interesting to complete    section 4.1 in
the setting of the iterated shtukas   and
Lafforgue's compactification for the shtuka  variety \cite{La1},
\cite{La2} and the arithmetic counterpart.

\bigskip

The paper is organized in the following way:

Section 2 is devoted to recalling some definitions and results of rank-$n$
vector bundles with level structures. In
subsection 2.2, we introduce  zeta correspondences for the moduli of
semistable  vector bundles with level structures
over ${\Bbb P}_1$, which are    analogous to the $n^2$-generalized theta
divisor. Moreover, zeta correspondences are
given by the zero locus of $n^2$-functions, obtained from a matrix.

In section 3. we recall some definitions and results for rank-$n$ Drinfeld
modules and elliptic sheaves with level
structures. Via the immersion of elliptic sheaves in  the moduli of
semistable  vector bundles, we gain a notion of zeta
correspondences for rank-$n$ Drinfeld modules. By means of a similar  Lemma
to Lemma 3.3.1 \cite{An2}, we can give some
information about the irreducible components of these zeta correspondences.
We provide rank-$2$ examples and the known
 rank-$1$ examples, \cite{An3}, \cite{AnDP}, \cite{Si}.

In section 4, by taking  a rational point $p$ instead of $\infty$, we can
obtain in a more precise way the results of the
last sections for a general curve $C$  and rank-$n$ shtukas. We state  a
similar Lemma    to Lemma 3.3.1 \cite{An2} for
rank-$n$
 shtukas, although in this case we should establish this Lemma   in rank
$n$ without
level structures to complete this Lemma.

\bigskip
\bigskip
 {\bf List of notations and previous definitions}
\bigskip

${{\Bbb F}_q}  $ is a finite field with $q$-elements, ($q=p^m$)

\bigskip
$\otimes$   denotes $ \underset  {{\Bbb F}_q} \otimes$

\bigskip
$R$ is a ${\Bbb F}_q$-algebra

\bigskip

 $C$       is  a smooth, proper and geometrically irreducible curve over
${\Bbb F}_q$

\bigskip
$g$          is the genus of $C$

\bigskip
$\infty $ and $p$ are rational points in $C$

\bigskip
$A$      $=H^0(C\setminus \{\infty  \},{{\cal O}_C})$

\bigskip
$D$ is an effective rational divisor over $Spec(A)$

\bigskip
If $M$ is a vector bundle over $C$, $M(k)$ denotes $M\underset {\o_C}
\otimes \o_C(k\infty)$, $k\in {\Bbb Z}$.

\bigskip

$R^\times$ denotes the group of units   in a ring $R$

\bigskip
If $s\in Spec(R)$, we denote by $m_s$ the prime ideal associated with $s$.

\bigskip

If $M$ is a vector bundle over $C\times Spec(R)$, and $s\in Spec(R)$,
$deg(M_s)$ denotes the degree of the vector bundle, $M_s$, over $C\times
k(s)$. In this paper
we shall consider  vector bundles of  constant degree, $deg(M)$, for each
$s \in Spec(R)$.

 \bigskip

Let  $M$ be  a vector bundle over $\o_C$;  $M^\vee$ denotes
the dual vector bundle

\bigskip

$({\Bbb G}_a)_R$      is   the additive line group over a ring  $R$
\bigskip

$F$ denotes the Frobenius morphism over a scheme $Spec(R)$

\bigskip
For a vector bundle, $M$, over $C\times Spec(R)$, $F^\#(M)$ denotes the
pull-back $(Id\times F)^*M$.

\bigskip
If $T:M\to M'$ is a morphism of vector bundles over a scheme, $H^0(T)$
denotes the morphism induced among the global
sections


\section{  Vector bundles
 with level structures    and zeta correspondences  over ${\Bbb P}_1$}

\subsection{ Vector bundles with level structures}

\begin{defn}
A $\infty  D$-level structure is a pair
$(M,f_{\infty  D})$, where $M$ is a rank-$n$ vector bundle  over $C\times
Spec(R)$  and
$$f_{\infty  D}:M\to (\o_C/\o_C(- \infty-D )\otimes R)^n $$
is a surjective morphism of $\o_C \otimes R$-modules.

By $ \o_C/\o_C(- \infty-D)    $ and $ O_{\infty D}$ we denote $ (\o_{
\infty D}) $
 and $H^0(C,(\o_C/\o_C(- \infty-D) ) $ respectively.

A morphism of two pairs, $(M,f_{\infty   D})$, $(M',f'_{\infty   D})$, is a
morphism of
$\o_C \otimes R$-modules, $\phi:M\to M'$, such that the diagram
$$ \xymatrix {    M   \ar[dr]^{\bar f_{\infty   D}}    \ar[r]^{\phi} &
  M'\ar[d]^{f'_{\infty   D}}
\\ & (\o_{- \infty-D})^n\otimes R } $$
is commutative.
Two pairs $(M,f_{\infty   D})$, $(M',f'_{\infty   D})$ are said to be
equivalent  if there exists an isomorphism $\phi$.
\end{defn}

\bigskip
\begin{propo}\label{mor} If $\phi:M\to M'$ is a morphism between two pairs
$(M,f_{\infty   D})$, $(M',f'_{\infty   D})$
 defined over a field $K$, then $\phi$ is injective.
\end{propo}
\begin{proof} Since   $M/Ker(\phi)\hookrightarrow M'$, $M/Ker(\phi)$ is a
free torsion module, and since
$$M/Ker(\phi)\underset {\o_C}\otimes \o_{\infty D}\simeq    (\o_{ \infty
D})^n \otimes K, $$
 $M/Ker(\phi)$ is
locally free of rank $n$, and hence $Ker(\phi)=0$
\end{proof}

We recall that $M$ is said to be semistable if for each geometric point
$s\in Spec(R)$, $M_s:=M\underset{R} \otimes
k(s)$ is semistable; i.e, for each $F\subset M_s$
$$\mu(F):=deg(F)/rank(F)\leq \mu(M_s)$$

Let denote us by ${\cal M}^{ss}_{C}(n,h,\infty D)$ the moduli scheme of
pairs $(M,f_{\infty   D})$  of semistable
vector bundles, $M$, of rank $n$, degree $h$, with a ${\infty   D}$-level
structure, over a curve $C$ (c.f :{\cite{Ss}}).
Then,
$$Hom_{Schem.}(Spec(R),{\cal M}^{ss}_{C}(n,h,\infty D))  =\{ \text{ pairs
over } R, (M,f_{\infty   D})
\} /
\text{Up to isomor.}$$

In the same way, we shall denote ${\cal M}_{C}(n,h,\infty D)$, the moduli
stack, {\cite {LM}}, of pairs
$(M,f_{\infty   D})$, of
 vector bundles, $M$, of rank $n$, degree $h$,   with a ${\infty
D}$-level structure.

\bigskip

In all this section, $C={\Bbb P}_1$. By choosing a local parameter $1/t$ in
$ \infty  $, $A={\Bbb F}_q[t]$.
We consider an isomorphism $H^0(\o_{{\Bbb P}_1},\o_{{\Bbb P}_1}/\o_{{\Bbb
P}_1}(-\infty - D   ))=
{\Bbb F}_q[t]/p(t)\times {\Bbb F}_q $, with
$p(t)=(t-\alpha_1)^{r_1}.(t-\alpha_2)^{r_2}\cdots (t-\alpha_l)^{r_l}$ and
$D=r_1.x_1+\cdots+r_l.x_l$, where $x_i$ is
the point associated with the maximal $(t-\alpha_i){\Bbb F}_q[t]$. Recall
that $d=deg(D)$.

\begin{propo} \label{sequence}Let $M$ be a semistable vector bundle over
${\Bbb P}_1\times Spec(R)$  of rank $n$
and degree
$0$ with a $\infty   D$-level structure $f_{\infty   D}$. Then:

1) $H^0({\Bbb P}_1\times Spec(R), M )$ is a free  $R$-module of rank $n$,
and $M=H^0({\Bbb P}_1\times Spec(R), M )\otimes
\o_C$.

2)Two pairs $(M,f_{\infty   D})$, $(\bar M,\bar f_{\infty   D})$,
$deg(M)=deg(\bar M)=0$, are equivalent  if and only
if by choosing bases
$\{s_1,\cdots,s_n\}$ and $\{\bar s_1,\cdots,\bar s_n\}$ for $H^0({\Bbb
P}_1\times Spec(R), M )$ and  $H^0({\Bbb
P}_1\times Spec(R), \bar M )$ respectively, there exists $g\in Gl_n(R)$, with
$$ f_{\infty   D}(s_i) = \bar f_{\infty   D}(g(\bar s_i))$$

\end{propo}

\begin{proof}
1) By taking global sections in the exact sequence of $\o_{{\Bbb P}_1}\otimes
R$-modules
$$0\to M(-\infty )\to M \to M/M(-\infty ) \to 0,$$
and bearing in mind the morphism given by the $D$-level structure $
f_\infty:M\to (\o_\infty  \otimes R)^n$, we obtain
an isomorphism
$$M/M(-\infty )\simeq (\o_\infty  \otimes R)^n$$
and we conclude because, as $M$ is semistable, we have
$$H^0({\Bbb P}_1\times Spec(R), M(-\infty ))=H^1({\Bbb P}_1\times Spec(R),
M(-\infty ))=0$$
 and therefore
 $H^0(f_\infty):H^0({\Bbb P}_1\times Spec(R), M )\to (O_\infty\otimes R)^n$
is an isomorphism.
The second assertion is easily deduced because $ H^0({\Bbb P}_1\times
Spec(R), M )\otimes
\o_C \subseteq M$ and $deg(M)=0$.
\bigskip

2) If $(M,f_{\infty   D})$, $(\bar M,\bar f_{\infty   D})$ are equivalent,
there exists an isomorphism of modules
$\phi:M\to \bar M$ such that $\bar f_{\infty   D}.\phi=f_{\infty   D}$.  In
this case,  $g\in Gl_n(R)$ is given by
considering
$H^0(\phi)$ and bases for $H^0({\Bbb P}_1\times Spec(R), M )$ and
$H^0({\Bbb P}_1\times Spec(R), \bar M )$.

Conversely,   if $\{s_1,\cdots,s_n\}$ and $\{\bar s_1,\cdots,\bar s_n\}$
are bases for
$H^0({\Bbb P}_1\times Spec(R),   M )$ and $H^0({\Bbb P}_1\times
Spec(R),\bar M )$ respectively, satisfying
$$ f_{\infty   D}(s_i) =  \bar f_{\infty   D}(g(\bar s_i))$$

for all
$i$, $\phi$ is obtained from $g$, since $g$ can be interpreted as an
isomorphism:
$$M\simeq (\o_{{\Bbb P}_1}\otimes R).s_1,\cdots,(\o_{{\Bbb P}_1}\otimes R).s_n
 \to \bar M \simeq (\o_{{\Bbb P}_1 }\otimes R).\bar s_1,\cdots,(\o_{{\Bbb
P}_1}\otimes R).\bar s_n$$

\end{proof}

\begin{rem}\label{universal}
 By considering a basis
$\{s_1,\cdots,s_n\}$, by the last Lemma $(M,f_{\infty   D})$ has
associated, in a bijective way, a matrix:
$$(\Delta_0+\Delta_1t+\cdots +\Delta_{d-1}t^{d-1})\times \Delta_\infty$$
where $\Delta_k, \Delta_\infty$ are $n\times n$-matrices   with entry
elements in $R$. These matrices are given by
$\{f_{\infty   D}(s_1),\cdots,f_{\infty   D}(s_n) \}\subset ( O_{\infty
D})^n\otimes R$.

One can consider the vector bundle (scheme) over ${\Bbb F}_q$,
$V(Hom_{{\Bbb F}_q}({\Bbb F}_q^n,  (O_{\infty D})^n))$
associated with the
 ${\Bbb F}_q$-vector space
$Hom_{{\Bbb F}_q}({\Bbb F}_q^n,  (O_{\infty D})^n)$. Then, for this vector
bundle there exists a  universal object
$$({\cal D}_0+{\cal D}_1t+\cdots +{\cal D}_{d-1}t^{d-1})\times {\cal D}_\infty$$
  ${\cal D}_k$ and  ${\cal D}_\infty$   being matrices with   the
independent variables $d_{i,j,k}$  and
$d_{i,j, \infty }$ as entries,
 respectively.
\end{rem}

\begin{lem}\label{functor} There exists an injective immersion (of schemes)
$$ \Psi:{\cal M}^{ss}_{{\Bbb P}_1}(n,0,\infty D) \longrightarrow
V(Hom_{{\Bbb F}_q}({\Bbb F}_q^n, (O_{\infty
D})^n))/Gl_n$$

where $Gl_n$  denotes the $n$-linear algebraic group  $/{\Bbb F}_q$ acting on
$$V(Hom_{{\Bbb F}_q}({\Bbb F}_q^n,   (O_{\infty D})^n))$$
  over the term ${\Bbb F}_q^n$.
\end{lem}
\begin{proof} Noticing that $V(Hom_{{\Bbb F}_q}({\Bbb F}_q^n, (O_{\infty
D})^n))/Gl_n $ is the Grassmannian of
$n$-planes in  $(O_{\infty   D})^n$, then $\Psi (M,f_{\infty   D})$ is
defined as the rank-$n$ subbundle of $ (O_{\infty
D})^n$ generated by
$\{f_{\infty   D}(s_1),\cdots,f_{\infty   D}(s_n)\}$.
If
$$(M,f_{\infty   D})\in {\cal M}^{ss}_{{\Bbb P}_1}(n,0,\infty D)(R),$$
then
 $<f_{\infty   D}(s_1),\cdots,f_{\infty   D}(s_n)>$ is a  rank-$n$
subbundle of $   (O_{\infty D})^n\otimes R$,
because taking sections in the exact sequence
$$0\to   M(-\infty-D)\to M \overset {f_{\infty   D}} \longrightarrow
(\o_\infty D)^n\otimes R\to 0$$
  one obtains the exact sequence
$$  <f_{\infty   D}(s_1),\cdots,f_{\infty   D}(s_n)> \hookrightarrow
(O_\infty D)^n\otimes R \to
H^1({\Bbb P}_1\times Spec(R),    M(-\infty-D)\to 0.$$
To obtain this last exact sequence we have used  the fact that $M$ is
semistable and hence $H^1({\Bbb P}_1\times
Spec(R),   M )=0$ and
$H^0({\Bbb P}_1\times Spec(R),   M(-\infty-D ))=0$. Moreover,  by
Grauert's theorem $H^1({\Bbb
P}_1\times Spec(R),    M(-\infty-D ))$ is a locally free
$R$-module.
\end{proof}

\begin{rem}
  $ \Psi$ takes values in  the open subscheme, $U$, of the vector bundle   of
$$Hom_{{\Bbb F}_q}({\Bbb F}_q^n,   (O_  D)^n)) $$
where $U$  are the morphisms   such that their composition with the natural
projection $(O_\infty D)^n \to (O_\infty  )^n$
is an isomophism. By fixing an isomophism ${\Bbb F}_q^n \to (O_\infty
)^n$, $U$ can be identified with the vector bundle
$$V(Hom_{{\Bbb F}_q}({\Bbb F}_q^n,   (O_  D)^n)).$$
Thus, by $ \Psi$, ${\cal M}^{ss}_{{\Bbb P}_1}(n,0,\infty D)$ is isomorphic
to $V(Hom_{{\Bbb F}_q}({\Bbb F}_q^n,   (O_
D)^n))$
\end{rem}

\bigskip


\subsection{Zeta correspondences  for vector bundles over ${\Bbb P}_1$ with
level structures}\label{pi}

In this section we shall construct   a   zeta correspondence for ${\cal
M}^{ss}_{{\Bbb P}_1}(n,0,\infty D)$. To do
so, we   take  into account the generalized theta function
\cite{An3} 6.1.1.

Over $\o_{{\Bbb P}_1}(k\infty)$ one can consider the following $\infty
D$-level structure, $\pi^k_{\infty   D}$:
$$\pi^k_{\infty   D}:=\pi^k_\infty\times \pi_D:\o_{{\Bbb P}_1}(k ) \to
\o_{{\Bbb P}_1}/\o_{{\Bbb
P}_1}(-1)\times \o_{{\Bbb P}_1}/\o_{{\Bbb P}_1}(-D)$$
 $\pi_D$ being  induced by the natural epimorphism $\o_{{\Bbb P}_1}\to
\o_{{\Bbb P}_1}/\o_{{\Bbb P}_1}(-D) $
and the natural inclusion $\o_{{\Bbb P}_1}\hookrightarrow
\o_{{\Bbb P}_1}(k )$. $\pi_\infty^k$ is obtained from the epimorphism
$$\o_{{\Bbb P}_1}(k )\to\o_{{\Bbb P}_1}(k )/\o_{{\Bbb P}_1}( k-1  )$$
and the isomorphism induced by the multiplication by $1/t^{k-1}$
$$\o_{{\Bbb P}_1}(k )/\o_{{\Bbb P}_1} ( k-1 )\simeq \o_{{\Bbb P}_1}
/\o_{{\Bbb P}_1}( -1).$$

In this way if $\lambda_0+\lambda_1.t+\cdots+\lambda_{k-1}.t^{k-1}\in
H^0({\Bbb P}_1,\o_{{\Bbb P}_1}(k ))$,
$$H^0(\pi_{\infty   D}^k)(\lambda_0+\lambda_1.t+\cdots+\lambda_{k-1}.t^{k-1})=
(\lambda_0+\lambda_1.t+\cdots+\lambda_{k-1}.t^{k-1},\lambda_{k-1})\in
O_D\times O_\infty$$

\begin{defn}\label{zeta} We define the $\infty  D$-zeta correspondence,
$Z^{\infty   D}_n$,
as the subscheme of   ${\cal M}^{ss}_{{\Bbb P}_1}(n,0,\infty D)\times {\cal
M}^{ss}_{{\Bbb P}_1}(n,0,\infty D)$ defined by

 $[(\bar M,\bar f_{\infty   D}),(M,f_{\infty   D})]\in
Hom_{schemes}(Spec(R),Z^{\infty   D}_n) $   if and only
if there exists a morphism of modules
$T:\bar M\to M(d-1) $ ($d=deg(D)$), where the diagram :

 $$ \xymatrix {  \bar M   \ar[dr]^{\bar f_{\infty   D}}    \ar[r]^{T} &
  M(d-1)\ar[d]^{f_{\infty   D} \otimes \pi^{d-1}_{\infty   D}}
\\ &(\o_\infty D)^n \otimes R } $$
is commutative.

Therefore, $[(\bar M,\bar f_{\infty   D}),(M,f_{\infty   D})]\in Z^{\infty
D}_n$ if and only if
$$Hom_{pairs}((\bar M, \bar f_{\infty   D}),(M(d-1 ),f_{\infty
D}\otimes\pi^{d-1}_{\infty   D}))\neq 0.$$

\end{defn}
From Lemma \ref{Zeta} it can be deduced  that this subscheme is closed.

\begin{lem} Given level structures $(\bar M,\bar f_{\infty
D}),(M,f_{\infty   D}) \in {\cal M}^{ss}_{{\Bbb
P}_1}(n,0,\infty D)$ over $R$, there exists a unique morphism of vector
bundles $T:\bar M\to M(d-1) $ with the diagram
$$ \xymatrix {  \bar M   \ar[dr]^{\bar f_{   D}}    \ar[r]^{T} &
  M(d-1)\ar[d]^{f_D \otimes \pi_D }
\\ &(\o_D)^n \otimes R } $$
 commutative.
\end{lem}
\begin{proof} By choosing bases
$\{s_1,\cdots,s_n\}$ and $\{\bar s_1,\cdots,\bar s_n\}$ for $H^0({\Bbb
P}_1\times Spec(R), M )$ and  $H^0({\Bbb
P}_1\times Spec(R), \bar M )$ respectively.

By   Remark \ref{universal}, $(M,f_{\infty   D})$, $(\bar M,\bar f_{\infty
D})$ has associated
$$(\Delta_0+\Delta_1t+\cdots +\Delta_{d-1}t^{d-1})\times \Delta_\infty$$
and
$$(\bar\Delta_0+\bar\Delta_1t+\cdots +\bar\Delta_{d-1}t^{d-1})\times
\bar\Delta_\infty,$$
respectively.

 Since
$$H^0({\Bbb P}_1\times Spec(R), M(d-1) )=\overset {d-1} {\underset{i=0}
\bigoplus} H^0({\Bbb P}_1\times Spec(R),
M ).t^i$$
by considering the above basis, $T$ is given by
$$  A_0+  A_1.t+\cdots +A_{d-1}.t^{d-1}, $$
where $A_i$ are $n\times n$-matrices with entries in $R$.

These matrices can be obtained from the relation, $\bar f_D= f_D.T_{\vert
D}$, where $T_{\vert D}$ is the isomorphism
$$T_{\vert D}:\bar M/\bar M(-D)\to M(d-1)/
M(d-1)(-D).$$
    Thus, $T$  is given by

$$  A_0+  A_1t+\cdots +A_{d-1}t^{d-1}=( \bar\Delta_0+ \bar \Delta_1t+\cdots
+ \bar\Delta_{d-1}t^{d-1})^{-1}.
( \Delta_0+ \Delta_1t+\cdots +\Delta_{d-1}t^{d-1}).$$

In this Lemma we have considered an isomorphism $O_D\simeq {\Bbb
F}_q[t]/p(t)$ and $A_i$ matrices for linear morphisms
between $H^0({\Bbb P}_1\times Spec(R),\bar M )$ and  $H^0({\Bbb P}_1\times
Spec(R),   M )$ expressed in the bases
$\{s_1,\cdots,s_n\}$ and $\{\bar s_1,\cdots,\bar s_n\}$.

\end{proof}

\begin{lem}\label{hom} With the above notation, $T:\bar M\to M(d-1) $ makes
the diagram
$$ \xymatrix {  \bar M   \ar[dr]^{\bar f_{\infty   D}}    \ar[r]^{T} &
  M(d-1)\ar[d]^{f_{\infty   D} \otimes \pi^{d-1}_{\infty   D}}
\\ &(\o_{\infty D})^n \otimes R } $$
  commutative  if and only if $\bar\Delta_\infty=\Delta_\infty.A_{d-1}$.
\end{lem}
\begin{proof} By following the same notation as in the latter lemma, let us
consider   $T=A_0+  A_1t+\cdots
+A_{d-1}t^{d-1}$. Then, the only condition that
$T$ must satisfy for   the diagram
$$ \xymatrix {  \bar M   \ar[dr]^{\bar f_{\infty   D}}    \ar[r]^{T} &
  M(d-1)\ar[d]^{f_{\infty   D} \otimes \pi^{d-1}_{\infty   D}}
\\ & (\o_\infty D)^n \otimes R } $$
to be commutative is that
$$ \xymatrix {  \bar M   \ar[dr]^{\bar f_\infty }    \ar[r]^{T} &
  M(d-1)\ar[d]^{{f_\infty } \otimes \pi^{d-1}_\infty }
\\ &(\o_\infty )^n\otimes R } $$
must be commutative. But  from the definition of $\pi^{d-1}_\infty$, if
$i<d-1$ then  $\pi^{d-1}_\infty(t^i)=0$  and hence
the condition that $A_0+  A_1t+\cdots +A_{d-1}t^{d-1}$ must be satisfied to
make  the above diagram commutative is
$\bar \Delta_\infty={\Delta_\infty}.A_{d-1}$.
\end{proof}

\begin{rem}\label{regular}
From Remark \ref{universal}
$$V(Hom_{{\Bbb F}_q}({\Bbb F}_q^n, (O_{\infty D})^n))=Spec({\Bbb
F}_q[\{d_{i,j,r}\}_{i,j,r},
\{w_{i,\infty  }\}_{i,\infty   }])$$
"$d_{i,j,r}$" and "$d_{i,j, \infty }$" being the matrix entries for ${\cal D}_r$
and
$ {\cal D}_\infty$, respectively, with
$$({\cal D}_0+{\cal D}_1t+\cdots +{\cal D}_{d-1}t^{d-1})\times {\cal D}_\infty$$
a universal object for $V:=V(Hom_{{\Bbb F}_q}({\Bbb F}_q^n, (O_\infty D)^n))$.
Let us denote  $p_1,p_2:V\times V\to V$ the natural projections and
${\cal A}_0+{\cal A}_1t+\cdots +{\cal A}_{d-1}t^{d-1}$ the matrix composition
 $$(p_2^*{\cal D}_0+p_2^*{\cal D}_1t+\cdots +p_2^*{\cal
D}_{d-1}t^{d-1})^{-1}.(p_1^*{\cal D}_0+\cdots +p_1^*{\cal D}_{d-1}t^{d-1}).$$
We now consider   $I^n_{\infty D}$, the subscheme of $V\times V$, given by
the relation
$$p_1^*{\cal D}_\infty=p_2^*{\cal D}_\infty.{\cal A}_{d-1}.$$
Since $I^n_{\infty D}$ is invariant, by the action of $Gl_n\times Gl_n $,
we obtain a subscheme
$$I^n_{\infty D}/Gl_n\times
Gl_n  \subset V/Gl_n\times V/Gl_n.$$
 Moreover, by the very definition of $I^n_{\infty D}$  it is given by the
zero locus of
$n^2$-regular functions of
$${\Bbb F}_q[\{d_{i,j,r}\}_{i,j,r}, \{d_{i,j,\infty  }\}_{i,j
}]\otimes{\Bbb F}_q[\{d_{i,j,r}\}_{i,j,r},
\{d_{i,\infty  }\}_{i,j   }]$$

\end{rem}

\bigskip

\begin{lem}\label{Zeta} Bearing in mind   the morphism $\Psi$ of
\ref{functor}, we have
$$(\Psi \times \Psi)^{-1} (I^n_{\infty D}/Gl_n\times Gl_n)=Z^{\infty   D}_n $$
and hence $Z^{\infty   D}_n $   is the zero locus of $n^2$-regular
functions of   ${\cal M}^{ss}_{{\Bbb P}_1}(n,0,\infty D) \times {\cal
M}^{ss}_{{\Bbb P}_1}(n,0,\infty D)$. Recall that the
subscheme
$Z^{\infty   D}_n$ is defined in \ref{zeta}.

\end{lem}
\begin{proof}
This Lemma is deduced from Lemma \ref{hom} and Remark \ref{regular}.

$Z^{\infty   D}_n $   is the zero locus of $n^2$-regular
functions of   ${\cal M}^{ss}_{{\Bbb P}_1}(n,0,\infty D) \times {\cal
M}^{ss}_{{\Bbb P}_1}(n,0,\infty D)$ because  of
Remark
\ref{regular}. Moreover, one can calculate in an explicit  way
these $n^2$-regular functions from the relation
$$0=(p_1^*{\cal D}_\infty)^{-1}.p_2^*{\cal D}_\infty.{\cal
A}_{d-1}$$
\end{proof}

\bigskip
For rank $1$, we have:

\begin{defn} We define the $\infty   D$-generalized theta divisor,
$\Theta_1^{\infty   D } $,
as the subscheme of   $J^{\infty   D}_{{\Bbb P}_1,0} $ defined by:
$(L,f_{\infty   D})\in Hom_{schemes}(Spec(R),\Theta_1^{\infty   D} ) $   if
and only if there exists a morphism of
modules
$T:\o_{{\Bbb P}_1} \otimes R\to L(d-1) $ ($d=deg(D)$), where the diagram:

$$ \xymatrix {  \o_{{\Bbb P}_1} \otimes R\ar[d]^{\pi_D }\ar[r]^{T} &
L(d-1)\ar[dl]^{f_{\infty   D}\otimes
\pi^{d-1}_{\infty   D} } \\
\o_{\infty   D}\otimes R & }$$ is commutative.

Here, $J^{\infty   D}_{{\Bbb P}_1,0} $ denotes the $\infty   D$-generalized
Jacobian  of line bundles of degree $0$,
\cite{Se}.
\end{defn}
Now, we   recall  the definition of the tensor product and dual level
structures for vector bundles.
\begin{rem} The tensor product of two $\infty   D$-level structures
$(M,f_{\infty   D})$ and
$(\bar M, \bar f_{\infty   D})$ is the $\infty   D$-level structure
$$(M\underset {\o_{{\Bbb P}_1}\otimes R}\otimes \bar M,f_{\infty
D}\otimes \bar f_{\infty   D})$$
where $f_{\infty   D}\otimes \bar f_{\infty   D}$ denotes the morphism
$$M\underset {\o_{{\Bbb P}_1}\otimes R}\otimes \bar M\to  (\o_{\infty   D} )^n
\underset{\o_{\infty   D}}
\otimes  (\o_{\infty   D} )^n  \otimes R\simeq  (\o_{\infty   D} )^{n^2}
\otimes R.$$

The dual $\infty   D$-level structure  $(M^\vee,f^\vee_{\infty   D})$   for
$(M,f_{\infty   D})$, is given by a level
structure, $f^\vee_{\infty   D}$, for the dual vector bundle $M^\vee$;i.e.,
  $f_{\infty   D}$ induces an isomorphism
$$M/M(-\infty - D)\simeq  (\o_{\infty   D} )^n   \otimes
R$$
 from which we deduce an isomorphism
$$M^\vee \underset {\o_{{\Bbb P}_1}\otimes R}\otimes  \o_{\infty   D}
\otimes R \to  (\o_{\infty   D} )^n   \otimes
R.$$
$f^\vee_{\infty   D}$ is given by this isomorphism.
\end{rem}

\bigskip
Let us now consider     the analogous morphism to the morphism \cite{An2}
4.1. but for line bundles with level
structures. Notice that the genus considered in the case of $\infty
D$-level structures is $g({\Bbb P}_1)+deg( \infty
D)-1=d$,
\cite{Se} pg.86.

Let $m_{\infty D}$ be the morphism
$$m_{\infty D}: J^{\infty   D}_{{\Bbb P}_1,0}  \times
 J^{\infty   D}_{{\Bbb P}_1,0}  \to  J^{\infty   D}_{{\Bbb P}_1,0} , $$
defined by
$$m_{\infty D}[(  L,  f_{\infty   D}),(\bar L,\bar f_{\infty   D})]=(\bar
L,\bar f_{\infty   D})\otimes( L^\vee ,
 f^\vee_{\infty   D}).$$
We have

\begin{lem}\label{m}
  For $n=1$, $m_{\infty D}^{-1}(\Theta_1^{\infty   D} )=Z^{\infty   D}_1$
and  $Z^{\infty   D}_1$ is a principal  Weil
divisor on
$J^{\infty   D}_{{\Bbb P}_1,0} \times    J^{\infty   D}_{{\Bbb P}_1,0}   $.
\end{lem}
\begin{proof} This is deduced from the latter definition, from the
definition of $Z^{\infty   D}_n$, and
from Lemma \ref{Zeta}.
\end{proof}

\section{Drinfeld modules and zeta correspondences}

\subsection{ $A$-Drinfeld modules, elliptic sheaves and their antiequivalence}

In this section we   recall the definition of Drinfeld modules, elliptic
sheaves, level structures
  and their antiequivalence \cite{ BlSt}, \cite{ Dr1}, \cite{
Dr2},\cite{LRSt}, \cite{L},\cite{Mu}.
\begin{defn} A Drinfeld module, $\phi$, of rank $n$ over $R$ is a ring morphism:
$$\phi:A\to End_R(({\Bbb G}_a)_R)$$
with $\phi_a=\underset {i=0} {\overset {nv_\infty(a)} \sum} a_i.\sigma^i$,
  $a\in A$, $v_\infty$ is the  valuation for $\infty$, $a_i \in R$,   $a_m$
is a unit in $R$ and  $({\Bbb
G}_a)_R$ is the additive line group over $R$, $End_R(({\Bbb G}_a)_R)$  are
its endomorphisms, and $\sigma$ is
the endomorphism on $({\Bbb G}_a)_R$,
$\sigma(\gamma)=\gamma^q$. In this way, $End_R(({\Bbb
G}_a)_R)=R\{\sigma\}$, where $R\{\sigma\}$ is the ring of
polynomials in $\sigma$, with the twisted rule of multiplication
$\sigma.b=b^q\sigma$.
\end{defn}

From this definition, one deduces  a morphism of rings
$$c:A\to R\{\sigma\} \overset {\sigma=0}\to R$$
called the characteristic morphism of $\phi$.

\begin{defn}\label{drin} An elliptic sheaf of rank $n$ over $R$,
$(E_{j},i_{j},\tau)$, is a commutative
diagram of vector bundles of rank
$n$ over $C\times Spec(R)$, and injective morphisms of modules $\{i_h\}_{ h\in
{\Bbb N}}$, $\tau$:

$$ \xymatrix @C=35pt{
E_1 \ar[r]^{i_1} &
E_2 \ar[r]^{i_2} &
 {\cdots} \ar[r]^{i_{(n-1)}}&
E_n  \ar[r]^{i_n} &
\cdots \\
 F^{\#} E_0\ar[u]^\tau \ar[r]^(.5){F^{\#} i_0} &
 F^{\#} E_1 \ar[u]^\tau \ar[r]^(.6){F^{\#} i_1} &
{\cdots }\ar[u]^\tau \ar[r]^(.6){F^{\#} i_n} &
 F^{\#} E_n \ar[u]^\tau \ar[r]& {\cdots}}$$

satisfying:

a) For any $z\in Spec(R)$, $deg((E_{h})_z)=ng+h$.
 \bigskip \bigskip

b) For all $i\in {\Bbb  Z}$, $E_{i+n}=E_{i}(1)$.
\bigskip \bigskip

c) $E_{i}+\tau(F^{\#} E_i )=E_{i+1}$.
\bigskip \bigskip

d) $j_*(E_i/E_{i-1})$ is a   rank-one free module over $R$. $j$ is the
inclusion $\infty \times Spec(R)
\hookrightarrow  C\times Spec(R)$.

 Recall that $F^{\#}$ denotes $(Id\times F)^*$.

 \end{defn}

\begin{rem}\label{trivial} From these properties, it is  deduced that
$h^0(E_{h})=h$ and $h^1(E_{h})=0$,
$h\geq 0$ c.f.{\cite{Dr2}}. In these conditions it is not hard to prove
that the "$E_{in}$" are semistable.

Moreover, it is proved   that there exists
    a basis $\{s,\sigma. s,\cdots, \sigma^{n+i-1}  . s\}$  for $H^0(C\times
Spec(R), E_i)$ ($i>0$), with
$\sigma. s:=\tau(F^{\#} s)$ and $\sigma^j. s:=\tau(F^{\#}\sigma^{j-1}.s)$.
\end{rem}

In {\cite{Dr2}},   the anti-equivalence between rank-$n$-Drinfeld modules
and rank-$n$-elliptic sheaves is settled
in the following way:
\begin{rem}\label{equivalence} Given a rank-$n$-Drinfeld module  over $R$
$$\phi:A\to End_R(({\Bbb
G}_a)_R)=R\{\sigma\}$$
$R\{\sigma\}$ has a structure of an
$A\otimes R$-module by means of $\phi$ ($R$-module on the left and
$A$-module on the right). Moreover, $R\{\sigma\}$
 has a graduation: $deg(\underset j  {\overset m  \sum} b_j.\sigma^j)=m$.
Bearing in mind this graduation,
the $A\otimes R$-module $R\{\sigma\}$ has associated a coherent sheaf $E_0$
over $C\times Spec(R)$. The "$E_j$"
are obtained  by translating the degree over $R\{\sigma\}$. $\tau$ is
obtained from the multiplication on the
left over $R\{\sigma\}$  by $\sigma$. It is not hard to prove that $E_i$
are vector bundles  of rank $n$.

Conversely, given an elliptic sheaf $(E_{j},i_{j},\tau)$ over $C\times
Spec(R)$, by regarding the diagram  associated with
this elliptic sheaf, for each
$h\in {\Bbb N}$ one deduces injective morphisms, which we denote by $\tau^h$
$$\tau^h: {F^\#}^{h }E_i\hookrightarrow E_{i+h} $$
and a chain of modules
$$\tau^k ({F^\#}^h E_i )\subset \tau^{k-1} ({F^\#}^{h-1}E_{i+1}
)\subset\cdots\subset
{F^\#}^{h-k} E_{h+k} $$ ($h\geq k$).

If  $<s>=H^0(C\times Spec(R), E_{1-n})$  by c) is
$$\underset i \bigcup H^0(C\times Spec(R), E_i)=R\{\sigma\}.s$$
with $\sigma^r.s =\tau^r.((F^{\#})^r(s))$,
then $a.s=\underset j {\overset m \sum} b_j.\sigma^j.s$, ($a\in A$) and
the Drinfeld  module associated with $(E_{j},i_{j},\tau)$ is
$\phi_a=\underset j {\overset m \sum} b_j.\sigma^j$.

Since $rank(E_i)=rank({F^\#}E_{i-1})$ and $deg(E_i)=deg({F^\#}E_{i-1})+1$,
$E_i/\tau({F^\#}E_{i-1} ) $ is a coherent
sheaf over $C\times Spec(R)$ such that for each $s\in Spec(R)$
$$E_i/\tau({F^\#}E_{i-1} ) \underset R \otimes k(s)$$
 is concentrated on $p_s\in Spec(A)$ and hence one obtains a morphism
$c^*:Spec(R)\to Spec(A)$. This morphism is the
characteristic morphism associated with $\phi$. We say that
$(E_{j},i_{j},\tau)$ has the generic characteristic if the
characteristic morphism $c:A\to R$ is injective.

\end{rem}

\begin{rem} Let $(E_{j},i_{j},\tau)$ be an elliptic sheaf defined over a
field $K$ and with generic characteristic
$c:A\hookrightarrow K$.
One can see that $ {F^\#}^{h }E_{i+1} /\tau({F^\#}^{h+1 }E_i) $
is concentrated on $y_h\in Spec(A)\times Spec(K)$,
with $y_0:=(Id\times F^h)(y_h)$ the characteristic of $(E_{j},i_{j},\tau)$
and $ E_n/\tau^n({F^\#}^{n}E_0 )$ is
concentrated on $y_0,\cdots,y_{n-1}\in Spec(A)\times Spec(K)$. Moreover, by considering the points "$y_i$" as morphisms
$$y_i:Spec(K)\to Spec(A)\times Spec(K) $$
they are given by the morphism $c^{q^i}:A\hookrightarrow K$ on the first
entry and the identity morphism on the second
entry. Therefore,
$y_i\neq y_j$ ($i\neq j$) because $c^{q^i}\neq  c^{q^j}$ since   $c$ is
injective.
\end{rem}

\begin{propo}\label{n} Let consider us two elliptic sheaves with generic
characteristic  $(E_{j},i_{j},\tau)$ and
$(\bar E_{j}, \bar i_{j},\bar \tau)$ defined over a field $K$. If $T:
E_0\to \bar E_{r}$ is a morphism of vector bundles
with the diagram
$$ \xymatrix {  E_0     \ar[rr]^{T} & &
 \bar E_{r}
\\{F^\#}^{n} E_{-n}\ar[u]^{\tau^n}\ar[rr]^{{F^\#}^{n}T (-1)}   &
&{F^\#}^{n}\bar E_{r-n}\ar[u]^{\bar \tau^n}   } $$
commutative, there exists a maximum $i\in {\Bbb N}$ such that
$T(E_0)\subseteq \bar \tau^i(
{F^\#}^{i}\bar E_{r-i})\subset \bar E_r $.  In this case, $(Id\times
F^i)y_0=z_0$, where $y_0$ and $z_0$ denote  the
characteristic of
$(E_{j},i_{j},\tau)$ and
$(\bar E_{j}, \bar i_{j},\bar \tau)$ respectively. We denote by $T(-1)$ the
morphism induced by $T$ over $E_{-n}:=E_0(-1)$.
\end{propo}
\begin{proof} If $i$ is the maximum $i\in {\Bbb N}$ such that
$T(E_0)\subseteq  \bar \tau^i( {F^\#}^{i}\bar E_{r-i})$, then we have a
non-trivial morphism
$$T: E_0/{F^\#}^{n}E_{-n}\to   \bar \tau^i( {F^\#}^{i}\bar E_{r-i})
/\bar\tau^{i+1}(  {F^\#}^{i+1}\bar E_{r-i-1})$$
where $\bar \tau^i( {F^\#}^{i}\bar E_{r-i}) /\bar\tau^{i+1}(
{F^\#}^{i+1}\bar E_{r-i-1})$ is concentrated on $z_i$ and
$E_0/{F^\#}^{n}E_{-n}$ is concentrated  on $y_0,\cdots,y_{n-1}$, $(Id\times
F^i)z_i=z_0$ and
$(Id\times F^j)y_j=y_0 $. We thus have
 $(Id\times F^j)z_0=(Id\times F^i)y_0$. To conclude it suffices to prove
that $j=0$.

We shall assume that $j\geq 1$.  However, in doing so we shall arrive at a
contradiction. It is not difficult to see that
each case
 is reduced to studying   $j=1$ and $i=0$. In this case, $z_0=y_1$ because
$(Id\times F )z_0=y_0=(Id\times F )y_1$. Thus,
$$T_{\vert {F^\#}E_{-1}} : {F^\#} E_{-1} \to    \bar E_r /\bar \tau( \bar
{F^\#} E_{r -1})$$
gives an isomorphism
$${F^\#}E_{-1}/\tau({F^\#}^2E_{-2} )\simeq    \bar E_r /\bar \tau(  {F^\#}
\bar E_{r -1})$$
(recall that ${F^\#}E_{-1}/\tau({F^\#}^{2}E_{-2} )$ is concentrated on
$y_1$). We therefore have a commutative diagram
$$ \xymatrix {  E_0  \ar[r]^{T}    \ar[r]^{T} &
 \bar E_{r}
\\ {F^\#}E_{-1}\ar@{^{(}->}[u]^\tau\ar[ur]_{{F^\#}T_{-1}}
&{F^\#}E_{r-1}\ar@{^{(}->}[u]^{\bar \tau }
\\ {F^\#}^2E_{-2} \ar@{^{(}->}[u]^\tau \ar[ur]  &   }, $$
where ${F^\#}T_{-1}$ is the restriction of $T$ to $ {F^\#}E_{-1}$. Since
$z_0=y_1$, we have  $z_{n-1}=y_{n }$ because
$ (Id\times F^{n-1})z_{n-1} =z_0=y_1=(Id\times F^{n-1})y_{n }$ and
$z_{n-1}$ and $y_n$ are given by injective morphisms
$A\hookrightarrow K$. Thus, in
the same way as before we also have  a commutative diagram
$$ \xymatrix {  &  {F^\#}^{n-1} \bar E_{r+1 -n} \\
{F^\#}^{n}E_{-n}  \ar[ur]^{{{F^\#}^nT}_{-n}} \ar[r]^{{F^\#}^{n}T(-1)}    &
{F^\#}^{n} \bar E_{r-n}
\ar@{^{(}->}[u]^{\bar
\tau}
 \\ {F^\#}^{n+1}E_{-n-1}\ar@{^{(}->}[u]^\tau \ar[ur] &  }, $$
 ${{F^\#}^nT}_{-n}$ being the restriction of $T$ to $ {F^\#}^nE_{-n}$.
Thus, this morphism gives an isomorphism
$${F^\#}^{n}E_{-n}/\tau( {F^\#}^{n+1} E_{-n-1} )\simeq     {F^\#}^{n-1}
\bar E_{r-n+1} /
\bar \tau( {F^\#}^{n }\bar E_{r-n} ).$$
However, this contradicts the inclusion
$${F^\#}^{n }T(-1)({F^\#}^{n }E_{-n})\subseteq   {F^\#}^{n } \bar E_{r-n} .$$
\end{proof}

\subsection{$A$-Drinfeld modules and elliptic sheaves with level structures}

Let  us now briefly recall    the definitions  of level structures for
elliptic sheaves and
Drinfeld modules:  $(E_{j},i_{j},\tau)$ and $\phi$  respectively.

Let $ E_{ I }$ be the subgroup scheme of $I$ division
points of $({\Bbb G}_a)_R$ as an $A$-module (via $\phi$).
$(\widetilde  {I^{-1}/A})^n$ will  denote the constant sheaf of stalk
$(I^{-1}/A)^n$.  $I$ is a proper ideal in $Spec(A)$.

\begin{defn}
An $I$-level structure for $\phi$  is a pair $(\phi ,\iota_{I})$.
$\iota_{I}$ is an  isomorphism of $A$-modules
$$\iota_{I}:E_{ I }(R)\to \widetilde {(  I^{-1}/A)^n}(R).$$

\end{defn}

\begin{defn}
A    $D$-level structure for the elliptic sheaf $(E_j,i_j,\tau)$,  is a
$D$-level structure, $f^j_D $, for each vector bundle $ E_j$
compatible with the morphisms
$\{i_j,\tau\}$. i.e., $ f^{D}_{rx} .i_j= f^j_D $ and
$f^{j+1}_D.\tau=(Id\times F)^*f^j_D$.
 \end{defn}

In {\cite{An1}},  the anti-equivalence  between Drinfeld modules and
elliptic sheaves with
level structures is established.  There exists an affine scheme, ${\cal
E}_n^D=Spec({\mathfrak R_n^D})$, that
represents rank-$n$ Elliptic sheaves with $D$-level structures
(equivalently, rank-$n$ Drinfeld  modules    with
$I$-level structures).

One can calculate in an explicit  way the sections $\sigma^i.s$ (Remark
\ref{trivial}) via   level structures. These
computations are done   in  {\cite{An1}}(Theorem 5), {\cite{Al}} (Remark 3.1).

Let $I$ be the ideal associated with the effective divisor $D\subset
Spec(A)$,   $D=r_1.x_1+\cdots+r_l.x_l$,
 $x_i$ being the point in $Spec(A)$ associated with the maximal ideal
$m_{x_i}$ and   $t_{x_i}$   a local uniformizer
for $m_{x_i}$.

\begin{rem}\label{calculation}  if $(\Phi,   \iota_I)$ is a  universal
Drinfeld module with an  $I$-level
structure  over ${\mathfrak R}_n^D$. Therefore, $  \iota_I$ is an isomorphism
$$  \iota_I:  E_I({\mathfrak R}_n^D)
\to (\underset {i\in \{1,\cdots,l\}} \prod   t_{x_i}^{-r_i}{\Bbb
F}_q[[t_{x_i}]]/  {\Bbb F}_q [[t_{x_i}]])^n$$
   with $  \iota_I(\alpha^{x_i}_{j h })=(0,\cdots,\overset{\overset j { \smile}}
{ v_j},\cdots,0)$
and $$v_j=(0,\cdots,\overset{\overset i { \smile}}
{ t^{-1-h }},\cdots,0)\in \underset {i\in \{1,\cdots,l\}} \prod
t_{x_i}^{-r_i}{\Bbb F}_q[[t_{x_i}]]/  {\Bbb
F}_q [[t_{x_i}]].$$

We denote by $(   {\Bbb E}_j,i_j,\tau,  f^j_D)$   the
 corresponding  universal elliptic sheaf of rank
$n$ with a  $D$-level structure associated with $(\Phi,   \iota_I)$. Thus,
  the value of the  sections $\sigma^k s \in H^0(C\otimes {\mathfrak
R}_n^D,     {\Bbb E}_m)$ via the level
structure  $  f_D^m:   {\Bbb E}_m\to \o_D^n\otimes {\mathfrak R}_n^D$ is:
$$ H^0( f_D^m)( \sigma^k s)=(s^k_1,\cdots,s^k_n)\in (O_D\otimes {\mathfrak
R}_n^D)^n$$
By using the isomorphism
$$(O_D\otimes {\mathfrak R}_n^D)^n=(\underset {i\in \{1,\cdots,l\}}
\prod   {\mathfrak R}_n^D[[t_{x_i}]]/ t_{x_i}^{ r_i} {\mathfrak R}_n^D
[[t_{x_i}]])^n$$
we have:
$$s^k_j=(\overset {r_1-1} {\underset {h  = 0}\sum}
  (\alpha^{x_1}_{j h})^{q^k} t_{x_1}^{h},\cdots,\overset {r_l-1} {\underset
{h= 0}\sum}
  (\alpha^{x_r}_{j h})^{q^k} t_{x_r}^{h}).$$
\end{rem}

\bigskip
\bigskip

\begin{defn}\label{infinity} A $\infty$-level structure for a rank-$n$
elliptic sheaf $(  E_j,i_j,\tau   )$ over
$R$  is a
$\infty$-level structure $ (   E_0, f_\infty)$ such that the diagram

$$ \xymatrix {  {F^\#}^{n}E_0    \ar[dr]^{  {F^{\#}}^n   f_\infty   }
\ar[r]^{\tau^n} &
 E_0(1)\ar[d]^{f_\infty   \otimes \pi_\infty }
\\ &(\o_\infty )^n \otimes R } $$
is commutative. $\pi_\infty$ is defined in \ref{pi}

To give a $\infty   D$-level structure for an elliptic sheaf is to give an
$\infty  $-level structure together with
a $  D$-level structure for the elliptic sheaf.
\end{defn}
The idea of considering $\infty$-level structures in this setting has been
suggested to me by G.W. Anderson.

As   for $D$-level structures, there exists an affine scheme, ${\cal
E}_n^{\infty D}=Spec({\mathfrak R_n^{\infty
D}})$, that represents rank-$n$ Elliptic sheaves with $ {\infty D}$-level
structures.

\bigskip

\begin{rem}\label{matrix}  Let $(  {\Bbb E}_j,i_j,\tau   )$ be a  universal
rank-$n$ elliptic sheaf with an $\infty
D$-level structure
$f_{\infty   D}=f_\infty \times f_D$ over $Spec({\mathfrak R}_n^{\infty
D})$. By choosing the basis
$$\{s,\sigma .s,\cdots,{\sigma}^{n-1} .  s\}$$
 in $H^0(C\otimes {\mathfrak R}_n^{\infty D}, {\Bbb E}_0)$,   Remark
\ref{universal} associates a matrix with  $({\Bbb E}_0, f_{\infty   D})$.
By usingthe results and notation of the latter
Remark,    this matrix is
$$(A_{x_1,0}+A_{x_1,1}.t_{x_1}+\cdots+A_{x_1,r_1-1}t_{x_1}^{r_1-1})\times
\cdots \times
(A_{x_l,0}+A_{x_l,1}.t_{x_l}+\cdots+A_{x_l,r_l-1}t_{x_l}^{r_l-1})\times
A_{\infty},$$

where $A_{x_i,h}$ is
$$\left(\begin{matrix}      \alpha^{x_i}_{1,h}  &(\alpha^{x_i}_{1,h})^q  &\cdots
&(\alpha^{x_i}_{1,h})^{q^{n-1}}
\\ \alpha^{x_i}_{2,h}  &(\alpha^{x_i}_{2,h})^q  &\cdots
&(\alpha^{x_i}_{2,h})^{q^{n-1}}
\\ \cdots &\cdots &\cdots &\cdots
\\  \alpha^{x_i}_{n,h}  &(\alpha^{x_i}_{n,h})^q  &\cdots
&(\alpha^{x_i}_{n,h})^{q^{n-1}}  \end{matrix}\right) .$$
$A_{\infty}$ is a $n\times n$-matrix with coefficients in  ${\mathfrak
R}_n^{\infty D}$ that
we shall calculate  in the case $C={\Bbb P}_1$ in the following subsection.
\end{rem}

\subsection{ ${\Bbb F}_q[t]$-Drinfeld  modules and zeta correspondences}
In the rest of this section   $C={\Bbb P}_1$.
We denote by ${\cal E}_n^{\infty   D}$ the moduli (scheme) of rank-$n$
elliptic sheaves with  $\infty   D$-level
structures for $A={\Bbb F}_q[t]$. We denote an object of this moduli by $(
E_j,i_j,\tau ,f_{\infty   D}   )$.

In a direct way, one can obtain a morphism
$$\Phi:{\cal E}_n ^{\infty   D}\to    {\cal M}^{ss}_{{\Bbb P}_1}(n,0,\infty D)$$
given by $\Phi(  E_j,i_j,\tau ,f_{\infty   D}   )=(  E_0,f_\infty \times
f^0_D   )$. Recall that by  Remark \ref{trivial},
$E_0$ are semistable.

  The next three subsections are devoted to studying  the zeta correspondences
 $$(\Phi\times \Phi)^{-1}(Z^{\infty   D}_n)$$
in
${\cal E}_n^{\infty   D}\times {\cal E}_n^{\infty   D}$, where $Z_{\infty
D}^n$ is defined as in  Definition
\ref{zeta}.

 $\Phi\times \Phi$ is the analogous morphism to that   defined in
\cite{An2} 4.1.
\begin{lem}\label{cuadrado} If
$$[(  E_j,i_j,\tau , f_{\infty D}),(  \bar E_j,\bar i_j,\bar \tau, \bar
f_{\infty D} )]\in (\Phi\times
\Phi)^{-1}(Z^{\infty   D}_n),$$
then there exists $h\in {\Bbb N}$ with $h\leq n(d-1)$ such that $(Id\times
F^h)y_0=z_0$, $y_0$ and $z_0$ being the
characteristics of
$(  E_j,i_j,\tau  )$ and $(  \bar E_j,\bar i_j,\bar \tau  )$ respectively.
\end{lem}
\begin{proof}
We assume the elliptic sheaves defined over a field $R$. By Definition
\ref{zeta},  there exists a diagram
(not necessarily commutative)
$$ \xymatrix {  E_0  \ar[rr]^T \ar[dr]      & &
 \bar E_{n(d-1)}\ar[dl]
\\ &(\o_{\infty D})^\otimes R &
\\ {F^\#}^{n}E_{-n}\ar[uu]^{\tau^n}\ar[rr]^{{F^\#}^{n}T(-1)} \ar[ur] &
&{F^\#}^{n}\bar E_{n(d-2)}\ar[uu]^{\bar
\tau^n}\ar[ul]    }
$$
(oblique arrows are given by level structures). This diagram  is
commutative after tensoring  by $\o_{\infty D }$
and we therefore have   that the morphism of vector bundles
$$T.\tau^n-\bar\tau^n.{F^\#}^{n}T(-1):{F^\#}^{n}E_{-n}\to  \bar E_{n(d-1)}$$
takes values in  $\bar E_{n(d-1)}(-\infty-D)$. However, we are in ${\Bbb
P}_1$ and $\bar E_{n(d-1)}$ and
${F^\#}^{n}E_{-n}$ are semistable. Therefore,
$$\bar E_{n(d-1)}\underset{{\Bbb P}_1\otimes R}\otimes
({F^\#}^{n}E_{-n})^\vee(-\infty-D)$$
is also semistable of degree $-n$. In this way,
$$H^0({\Bbb P}_1\otimes R,\bar E_{n(d-1)}\underset{{\Bbb P}_1\otimes
R}\otimes ({F^\#}^{n}E_{-n})^\vee(-\infty-D))=0.$$
Hence $T.\tau^n-\bar\tau^n.{F^\#}^{n}T(-1)=0$.
We therefore   conclude because we are in the conditions of   Proposition
\ref{n}.

\end{proof}

Let $K$ be the ${\Bbb F}_q(t)$-algebra, ($t\to x$), $K:={\Bbb
F}_q(x)[a_1,\cdots,a_{n-1}]$. Let $\phi$ be the rank-$n$
Drinfeld  module with generic characteristic, defined over $K$:
$$\phi_t=\sigma^n+a_{n-1}.\sigma^{n-1}+\cdots+a_1.\sigma+x.$$

Let $(  E_j,i_j,\tau )$ be an   elliptic sheaf associated with $\phi$ and
$s\in H^0({\Bbb P}_1\otimes K,E_{1-n})$, with
$t.s=x.s+a_1.\sigma .s+\cdots+a_{n-1}.\sigma^{n-1}.s+
\sigma^n .s$.

Bearing in mind   Definition \ref{infinity}, we   study  when an
$\infty$-level structure $f_\infty$ for
$E_0$ is an $\infty$-level structure for the  elliptic sheaf $(
E_j,i_j,\tau     )$.

\bigskip
Since $H^0({\Bbb P}_1\otimes K, E_0(1))=H^0({\Bbb P}_1\otimes K, E_0 )\oplus t.H^0({\Bbb P}_1\otimes K, E_0 )$
 if we consider the bases $\{s,\sigma .s,\cdots,{\sigma}^{n-1} .  s\}$ and
$\{{F^{\#}}^n(s),\sigma
.{F^{\#}}^n(s),\cdots, {\sigma}^{n-1} .  {F^{\#}}^n(s)\}$ for $H^0({\Bbb
P}_1\otimes K, E_0 )$ and $H^0({\Bbb
P}_1\otimes K, F^{\#} E_0 )$ respectively, the morphism
$\tau^n:{F^{\#}}^nE_0\to E_0(1)$ is given in these basis by
  $A.t+B$, $A,B$ being matrices with entries in $K$.

\begin{lem} We have that
$$A.t+B={F^{\#}}^{n-1}C.{F^{\#}}^{n-2}C...{F^{\#}}C.C$$
with
$$C=\left(\begin{matrix}   0      & 0 &\cdots & 0& t-x
\\   1 & 0 &\cdots &0 &-a_{n-1}
\\ 0&1 &\cdots &0 &-a_{n-2}
\\ \cdots   &\cdots &\cdots  &\cdots &
\\   0 &0 & \cdots & 1&-a_1  \end{matrix}\right) $$
\end{lem}

\begin{proof} Bearing in mind Remark \ref{equivalence},   $\tau^n$ is given by
$${F^*}^nK\{\sigma\}\overset{\sigma\overset n {...} \sigma} \longrightarrow
K\{\sigma\},$$
where $\sigma.$ is the multiplication on the left over $K\{\sigma\}$. In
this setting, the bases
$\{s,\sigma .s,\cdots,{\sigma}^{n-1} .  s\}$ and $\{{F^{\#}}^n(s),
{F^{\#}}^n(\sigma.s),\cdots,   {F^{\#}}^n({\sigma}^{n-1} .s)\}$ are
$\{1,\sigma ,\cdots,{\sigma}^{n-1} \}$ and
$\{{F^{\#}}^n(1),
{F^{\#}}^n(\sigma ),\cdots,   {F^{\#}}^n({\sigma}^{n-1}  )\}$ respectively.
We conclude bearing in mind that
$t.s=x.s+a_1.\sigma .s+\cdots+a_{n-1}.\sigma^{n-1}.s+
\sigma^n .s$.
\end{proof}

Now we shall calculate $A_\infty$ (Remark \ref{matrix})
\begin{lem}\label{A} Let $(  E_j,i_j,\tau )$ be a  rank-$n$ elliptic sheaf,
$ f_\infty $   an $\infty$-level
structure over $ E_0$, and $(\nu_{ij})_{i,j}$ the matrix associated with
$$f_\infty:H^0({\Bbb P}_1\otimes K, E_0 )\to K^{\oplus n}$$
for the basis $\{s,\sigma .s,\cdots,{\sigma}^{n-1} .  s\}$ and the standard
basis in $K^n$. Then $(E_0,f_\infty)$
is an $\infty$-level structure for the elliptic sheaf $(  E_j,i_j,\tau )$
if and only if
$(\nu_{ij}^{q^n} )_{i,j}= (\nu_{ij})_{i,j}.A$, where $A$ is defined as in
the above Lemma.
\end{lem}
\begin{proof}Because of the definition of $\infty$-level structures for
elliptic sheaves,  the
diagram
$$ \xymatrix {  {F^\#}^{n}E_0    \ar[dr]^{  {F{\#}}^n   f_\infty   }
\ar[r]^{\tau^n} &
 E_0(1)\ar[d]^{f_\infty   \otimes \pi_\infty }
\\ &(\o_\infty )^n \otimes R } $$
must be commutative.

Bearing in mind  1) of Proposition \ref{sequence} and the above Lemma, we
conclude by taking the bases
 $\{{F^{\#}}^n(s),
{F^{\#}}^n(\sigma.s),\cdots,   {F^{\#}}^n({\sigma}^{n-1} .s)\}$ in
$H^0({\Bbb P}_1\otimes K, {F^\#}^{n}E_0 )$ and
$$\{s,\sigma
.s,\cdots,{\sigma}^{n-1} .  s\}\cup t.\{s,\sigma .s,\cdots,{\sigma}^{n-1} .
s\}$$
for
 $H^0({\Bbb P}_1\otimes
K, E_0(1) )$.
\end{proof}
$A_\infty$ (Remark \ref{matrix}) is $(\nu_{ij})_{i,j}$.

\bigskip
We shall now consider $\bar K_n^{\infty   D}:=K[\{\alpha^{x_i}_{j
h}\}_{\{i,j,h\}}, \{\nu_{ij}\}_{i,j}]$ as a
subring of an algebraic closed field containing $K$, where the
"$\alpha^{x_i}_{j h}$" are defined as in   Remark
\ref{calculation} and "$\nu_{ij}
$"   are solutions for the equation $(\nu_{ij}^{q^n} )_{i,j}=
(\nu_{ij})_{i,j}.A$.

It is not hard to prove that there exists a morphism $Spec(\bar K_n^{\infty
D})\to {\cal E}_n^{\infty   D}$.

\begin{thm}\label{comp} If we denote by $\beta$   the morphism composition
$$ Spec(\bar K_{n,\infty   D})\times Spec(\bar K_{n,\infty   D})\longrightarrow
{\cal E}_n^{\infty   D}\times {\cal
E}_n^{\infty   D}\overset {\Phi \times \Phi}\longrightarrow {\cal
M}^{ss}_{{\Bbb P}_1}(n,0,\infty D)\times {\cal
M}^{ss}_{{\Bbb P}_1}(n,0,\infty D),$$
  then
$$\beta^{-1}(Z^{\infty   D}_n)$$
is the zero locus of the $n^2$-functions deduced from an equation
$$ \bar \Gamma_\infty. \bar\Gamma_{d-1}  -Id_{n }=0,$$
 $\bar\Gamma_{d-1}, \bar \Gamma_\infty  $ being   $n\times n$-matrices with
entries in
$ \bar K_n^{\infty   D} \otimes \bar K_n^{\infty   D} $.
\end{thm}
\begin{proof}This Theorem   is a direct consequence of Lemma   \ref{Zeta}
and one can make an   explicit  calculation
of $\bar\Gamma_{d-1}, \bar \Gamma_\infty
$  bearing in mind   Remark \ref{matrix} for  parameters $t_{x_i}=t-\alpha_i$, Lemma
  \ref{A}, and  the  Chinese
remainder theorem: i.e., there exists a ring isomorphism
$$\delta:     \underset {i\in \{1,\cdots,l\}}
\prod   {\Bbb F}_q[t ]/ (t-\alpha_i)^{ r_i} \to O_D:={\Bbb F}_q[t]/(p(t))$$
given by
$$\delta(h_1(t),\cdots,h_l(t))=
p_1(t-\alpha_1).h_1(t).p(t)/(t-\alpha_1)^{
r_1}+\cdots+p_l(t-\alpha_l).h_l(t).p(t)/(t-\alpha_l)^{ r_l} $$
 with
$deg (p_i(t-\alpha_i))< r_i$  and
$$1/p(t)=\underset {i\in \{1,\cdots,l\}} \sum
p_i(t-\alpha_i)/(t-\alpha_i)^{ r_i}.$$
Recall that $D=(p(t))_0$, with $p(t)=\underset {i\in \{1,\cdots,l\}}
\prod   (t-\alpha_i)^{ r_i}$.

\end{proof}

 \subsection{Some explicit  irreducible components for   zeta correspondences}

The finite group $Gl_{\o_{\infty   D}}(\o_{\infty   D}^n)=Gl_n({\Bbb
F}_q)_\infty\times Gl_{\o_D}(\o_D^n)$ acts on
$\bar K_n^{\infty   D}$ by acting on the level structures. By considering
the restrictions over $\bar K_n^{\infty   D}$
 $$a_1=q_1(x),\cdots,a_{n-1}=q_{n-1}(x)$$
  with $q_1(x),\cdots, q_{n-1}(x) \in {\Bbb F}_q[x]$, if we   denote
$$\tilde K_n^{\infty   D}:=\bar K_n^{\infty
D}/(a_1-q_1(x),\cdots,a_{n-1}-q_{n-1}(x))$$
we have
$$Spec(\tilde K_n^{\infty   D})/Gl_n({\Bbb F}_q)_\infty\times
Gl_{\o_D}(\o_D^n)=Spec({\Bbb F}_q(x)).$$

\begin{lem}\label{beta} Let $\tilde \beta$ be the restriction of $  \beta $
to $Spec(\tilde K_n^{\infty   D})\times
Spec(\tilde K_n^{\infty   D})$:
$$  \tilde \beta :Spec(   \tilde K_n^{\infty   D})\times Spec( \tilde
K_n^{\infty   D})\longrightarrow {\cal
M}^{ss}_{{\Bbb P}_1}(n,0,\infty D)
\times {\cal M}^{ss}_{{\Bbb P}_1}(n,0,\infty D).$$
Then,

$$( \tilde \beta )^{-1}(Z^{\infty   D}_n  )=\underset{\underset { \text{
and } 0\leq i \leq n(d-1)}{\text{ For certain }
g\in Gl_{\o_{\infty   D}}(\o_{\infty   D}^n)}} \bigcup \Gamma_{g.F^i}.$$
$F$ is the Frobenius morphism and $\Gamma_{g.F^i}$ is the graph for the
morphism $g.F^i$. $g$ is considered as an
automorphism over $\tilde K_n^{\infty   D}$.
\end{lem}
\begin{proof}
 Let $(  E_j,i_j,\tau, f_{\infty   D})$ and $(  \bar E_j,\bar i_j,\bar\tau,
\bar f_{\infty   D} )$ be two rank-$n$
elliptic sheaves with an $\infty   D$-level structure defined over $ \tilde
K_n^{\infty   D}$.
For Proposition \ref{n}, if
$$[(E_0,f_{\infty   D}),(\bar E_0,\bar f_{\infty   D})] \in ( \tilde \beta
)^{-1}(Z^{\infty   D}_n ),$$
then $F^j$ (Characteristic$(  E_j,i_j,\tau, f_{\infty   D}))=$
Characteristic$(  \bar E_j,\bar i_j,\bar\tau, \bar
f_{\infty   D} )$ for some   $j\in {\Bbb N}$. Thus,   for some $g\in
Gl_{\o_{\infty   D}}(\o_{\infty   D}^n)$  is
$$g.F^j(  E_j,i_j,\tau, f_{\infty   D}))= (  \bar E_j,\bar i_j,\bar\tau, \bar
f_{\infty   D} )$$
 because
$$Spec(  \tilde K_n^{\infty   D})\to Spec(\tilde K_n^{\infty
D})/Gl_{\o_{\infty   D}}(\o_{\infty   D}^n)=Spec({\Bbb
F}_q(x))$$
is the characteristic morphism.

\end{proof}

\begin{lem}\label{irreducible} Let $s(t)$ be a monic polynomial with
$d-1-i=deg(s(t))\leq d-1$. We denote by
$h_{s(t)}$ the central matrix
$Id_n\times diag(s(t),\cdots, s(t)) \in Gl_n({\Bbb F}_q)_\infty\times
Gl_{\o_D}(\o_D^n)$. Then,
$$\Gamma_{h_{s(t)}.F^{ni}}\subset (  \tilde \beta )^{-1}(Z^{\infty   D}_n).$$
Moreover, if
$$\Gamma_{g.F^{n(d-1)}}\subset ( \tilde\beta )^{-1}(Z^{\infty   D}_n),$$
then $g=Id$.

\end{lem}
\begin{proof}
It suffices to find $T:{F^{\#}}^{ni}(E_0)\to E_0(d-1)$ such that the diagram
$$ \xymatrix {  {F^{\#}}^{ni}(E_0)  \ar[dr]_{h_{s(t)}.{F^{\#}}^{ni}
f_{\infty   D}  }  \ar[r]^{T} &
  E_0(d-1)\ar[d]^{f_{\infty   D} \otimes \pi^{d-1}_{\infty   D}}
\\ &(\o_{\infty   D})^n \otimes   \tilde K_{n }^{\infty   D}  } $$
is commutative. However, since
$$ \xymatrix {   {F^{\#}}^{ni}(E_0 )   \ar[dr]_{{F^{\#}}^{ni} f_{\infty
D} \otimes
\pi^i_{\infty   D}}
\ar[r]^{\tau^{ni}} &
  E_0(i)\ar[d]^{f_{\infty   D} \otimes \pi^{d-1}_{\infty   D}}
\\ &(\o_{\infty   D})^n \otimes   \tilde K_{n}^{\infty   D}  } $$
 is commutative because $f_{\infty   D}$ is an $\infty   D$-level structure for the elliptic sheaf
$(  E_j,i_j,\tau)$, it suffices to find $T': E_0(i)\to E_0(d-1)$ with
$$ \xymatrix {   E_0(i)   \ar[dr]_{ h_{s(t)}.(f_{\infty   D} )\otimes
\pi^i_{\infty   D}}
\ar[r]^{T'} &
  E_0(d-1)\ar[d]^{f_{\infty   D} \otimes \pi^{d-1}_{\infty   D}}
\\ &(\o_{\infty   D})^n \otimes \tilde  K_{n}^{\infty   D}  } $$
commutative.  We shall take  $T'$   defined as the homothety by $s(t)$. It
is clear that it is defined from
$E_0(i)$ to
$E_0(d-1)$, because $deg(s(t))=d-1-i$. Moreover, the above diagram is
commutative because
$$ \xymatrix {   E_0(i)   \ar[dr]_{h_{s(t)}.( f_D )\otimes
\pi_D}
\ar[r]^{T'} &
   E_0(d-1)\ar[d]^{f_{    D} \otimes \pi_{    D}}
\\ &(\o_{    D})^n \otimes  \tilde K_{n}^{\infty   D}  } $$
is commutative by the definition of  $T'$. And
$$ \xymatrix {   E_0(i)     \ar[dr]_{h_{s(t)}.(f_\infty ) \otimes
\pi^i_\infty    }
\ar[r]^{T'} &
   E_0(d-1)\ar[d]^{f_{\infty    } \otimes \pi^{d-1}_{\infty    }}
\\ &(\o_{\infty    })^n \otimes  \tilde K_{n}^{\infty  D }  } $$
is commutative because  $s(t)$ is monic.

The second assertion of the Lemma is deduced  because $T$ takes values in
$$\tau^{n(d-1)}(E^{\sigma^{n(d-1)}}_0  )\subset E_{n(d-1)}$$
 because of Proposition \ref{n}.
  We thus have $\tau^{n(d-1)}=T$,  (up to isomorphisms). In this case, if
$$(f_{\infty   D} \otimes \pi^{d-1}_{\infty
D}).\tau^{n(d-1)}=g.{F^{\#}}^{n(d-1)} f_{\infty   D}, $$
this implies that $g=Id$ because   the definition of $\infty   D$-level
structure for an elliptic sheaf is
$$(f_{\infty   D} \otimes \pi^{d-1}_{\infty   D}).\tau^{n(d-1)}=
{F^{\#}}^{n(d-1)} f_{\infty   D} .$$

\end{proof}

From this Lemma we deduce:

\begin{thm} (c.f \cite{An3},\cite{AnDP}, \cite{Si})\label{n=1} For $n=1$
$$( \tilde \beta_1 )^{-1}(\Theta^{\infty   D}_1)=\underset{\underset {
0\leq i \leq d-1}
 {\underset{  s(t) \text{ monic with }   deg(s(t))=i }  {s(t)\in ({\Bbb
F}_q[t]/p(t))^\times }}}\bigcup
\Gamma_{s(t).F^i}$$
\end{thm}
Here, ${\tilde\beta}_1$ denotes $m_{\infty D}.\tilde\beta$ and we follow
the notation of Lemma \ref{m}.
\bigskip

We can now obtain the same result as in \cite{An3} merely by considering
all possibilities over the $\infty$-level
structures: From Lemma \ref{A} for rank $1$, two $\infty$-level structures
for a rank-$1$ elliptic sheaf differ  in
a
$a\in {\Bbb F}_q^\times$. Thus,
$$\underset {a\in {\Bbb F}_q^\times}\bigcup ( \tilde\beta_1)^{-1}(
\Theta_{1, a\times 1}^{\infty   D})=
\underset{\underset {   0\leq i \leq d-1, deg(s(t))=i}
 {s(t)\in ({\Bbb F}_q[t]/p(t))^\times }}\bigcup
\Gamma_{s(t).F^i}$$
Here, $a\times 1 \in {\Bbb F}_q^\times \times \o_D^\times$, $a\times 1$
acts on
${\cal M}^{ss}_{{\Bbb P}_1}(1,0,\infty D)$ (this scheme is the generalized
Jacobian   $J^{\infty   D}_{{\Bbb P}_1}$), and
$\Theta_{1, a\times 1}^{\infty   D}$ denotes the transformation of
$\Theta_1^{\infty   D}$ by   $a\times 1$.

\bigskip
\bigskip

\subsection{Examples}
As above,   $C$ will be ${\Bbb P}_1$.

\bigskip

{\bf Example. 1} (\cite{An4},\cite{C}) We shall study  $n=1$,
$p(t)=t(t-1)$. We thus denote $D$ by $0+1$,  and
$\tilde K_1^{\infty +0+1}={\Bbb F}_q(x)[ \alpha_0, \alpha_1, \nu]$,  where
$$\alpha_0^q+\alpha_0.x=0 \text{, }\alpha_1^q+\alpha_1.(x-1)=0,$$
and $\nu \in {\Bbb F}_q^\times$. We can set $\nu=1$. With the variable
changes $u=\alpha_0$ and $v=\alpha_1$, we have
$${\Bbb F}_q(x)[ \alpha_0, \alpha_1, \nu]={\Bbb F}_q(u)[v ,1/v
]/u^{q-1}-v^{q-1}+1$$

Now bearing in mind   the  Chinese remainder theorem, there exists a ring
isomorphism
$$\delta:
 {\Bbb F}_q[t ]/ ( t ) \times {\Bbb F}_q[t ]/  (t-1)  ={\Bbb F}_q[t]/(t(t-1)),$$
given by
$$\delta(c,d)=ct+d(1-t).$$
 In this way, $ct+d(1-t)$ is monic $\in ({\Bbb F}_q[t]/(t(t-1)))^\times$ if
and only if $c-d=1$ and $c\neq 0,1$.

For Theorem \ref{n=1},
$$( \tilde\beta_1)^{-1}(\Theta^{\infty +0+1}_1)=
\underset{\underset{ c\neq 0,1}{  c\in {\Bbb F}_q^\times}}
\bigcup \Gamma_{[c,c-1] } \cup \Gamma_F,$$
 $\Gamma_{[c,c-1]}$ being the graph of the morphism in $\tilde K_{1}^{\infty +0+1}$  defined by $u\to c.u$, $ v\to
(c-1).v$. and
$$ \tilde\beta_1:Spec(  \tilde K_{1}^{\infty +0+1})\times Spec( \tilde
K_{1}^{\infty +0+1})\longrightarrow
{\cal M}^{ss}_{{\Bbb P}_1}(1,0,\infty +0+1)=J_{\infty +0+1}^{{\Bbb P}_1},$$
defined in \ref{n=1}.
By   Theorem \ref{comp} and by the explicit  calculations in   Remark
\ref{matrix}  and Lemma
\ref{A}, $\tilde \beta_1 $ has associated the $1\times
1$-"matrix"
$$(1\otimes u\times 1\otimes v\times 1\otimes1)^{-1}.(u\otimes 1\times
v\otimes 1\times 1\otimes1) $$
and
$$( \tilde \beta_1 )^{-1}(\Theta^{\infty +0+1}_1) =(\frac{u\otimes
1}{1\otimes u}-
\frac{v\otimes 1}{1\otimes v}-1\otimes 1)_0.$$

If we consider $\nu =a\in {\Bbb F}_q^\times$, we obtain
$$( \tilde\beta_1)^{-1}(\Theta^{\infty +0+1}_{1,1\times a})
=(\frac{u\otimes 1}{1\otimes u}-
\frac{v\otimes 1}{1\otimes v}-a\otimes 1)_0 $$
and
$$( \tilde\beta_1)^{-1}(\Theta_{\infty +0+1}^{1,1\times a})=
\underset{\underset{ c\neq 0,a}{  c\in {\Bbb F}_q^\times}}
\bigcup \Gamma_{[c,c-a] } \cup \Gamma_F$$

\bigskip

{\bf Example. 2} In this example we shall study  the last example in a more
general way.    $n=1$,
$p(t)=\underset {i\in \{1,\cdots,l\}}\prod    (t-\alpha_i)$. Then,
$D=\alpha_1+\cdots+\alpha_l$ and
$$\tilde K_{1}^{\infty +D}={\Bbb F}_q(x) [\alpha_1,\cdots \alpha_l, \nu],$$
 where $\alpha_i^q+\alpha_i.(x-\alpha_i)=0 $ with $1\leq i \leq l$ and $\nu
\in {\Bbb F}_q^\times$.  We can set
$\nu=1$ as in the above example.

Let us consider  the decomposition in simple fractions
$$1/p(t)=\underset {i\in \{1,\cdots,l\}} \sum m_i/ t-\alpha_i ,$$
$m_i\in {\Bbb F}_q^\times$. The  Chinese remainder theorem is then settled
in the following way: There exists a
ring isomorphism
$$\delta:     \underset {i\in \{1,\cdots,l\}}\prod   {\Bbb F}_q[t ]/
(t-\alpha_i) \to ={\Bbb F}_q[t]/(p(t))$$
given by
$$\delta(h_1 ,\cdots,h_l )= m_1.h_1 .p(t)/ t-\alpha_1 +\cdots+ m_l .h_l
.p(t)/ t-\alpha_l  .$$
In this way, $\delta(h_1 ,\cdots,h_l ) \in ({\Bbb F}_q[t]/(p(t)))^\times$ if
  $h_i\neq 0$ for all $i$.

 For Theorem \ref{n=1}
$( \tilde\beta_1)^{-1}(\Theta^{\infty +\alpha_1+\cdots+\alpha_l}_1)$
is
$$\underset { 0\leq i \leq deg(p(t) }\bigcup
\underset{\underset {\delta(h_1 ,\cdots,h_l )\text{monic and }
h_1...h_l\neq 0}{\underset{ deg(\delta(h_1 ,\cdots,h_l
))=deg(p(t)-i }{  (h_1 ,\cdots,h_l )}}}
\bigcup \Gamma_{\delta(h_1 ,\cdots,h_l ).F^i }  $$

$\Gamma_{\delta(h_1 ,\cdots,h_l ).F^i}$ is the graph  of the morphism in
$\tilde K_{1}^{\infty +0+1}$  defined by
$\alpha_r\to h_r.\alpha_r^{q^i}$, for each $i$.

By   Theorem \ref{comp} and by the explicit  calculations in  Remark
\ref{matrix}  and Lemma
\ref{A}, $ \tilde\beta_1$ has
associated the $1\times 1$-"matrix"
$$(1\otimes \alpha_1\times\cdots \times
1\otimes \alpha_l\times 1\otimes1)^{-1}.(\alpha_1\otimes 1\times\cdots
\times \alpha_l\otimes 1\times 1\otimes1)$$
and
$$(\tilde \beta_1)^{-1}(\Theta^{\infty +\alpha_1+\cdots+\alpha_l}_1)
=(m_1.\frac{\alpha_1\otimes 1}{1\otimes
\alpha_1}+\cdots+
m_l.\frac{\alpha_l\otimes 1}{1\otimes \alpha_l}-1\otimes 1)_0.$$

\bigskip

{\bf Example. 3} We now consider     $n=2$,
$p(t)=t $. Then, $D=0$ and we follow  the notation of  Remark \ref{matrix}
$$\tilde K_{2}^{\infty +0 }={\Bbb F}_q(x)[ \alpha^0_{1,0}, \alpha^0_{2,0},
\nu_{11},\nu_{21},\nu_{12},\nu_{22}]$$
  where
$$(\alpha^0_{i,0})^{q^2} +q_1(x).(\alpha^0_{i,0})^q + \alpha^0_{i,0}.x=0 $$
for $i=1,2$ and $(\nu_{ij})_{ij}$ satisfies
$(\nu^{q^2}_{ij})_{ij}=(\nu_{ij})_{ij}.A$,
with $A$ defined as in   Lemma \ref{A}
$$At+B=\left(\begin{matrix}       0& t-x^q
\\   1 &  -q_1(x)^q
  \end{matrix}\right).\left(\begin{matrix}       0& t-x
\\   1 &  -q_1(x)
  \end{matrix}\right). $$
Thus,
$$A=\left(\begin{matrix}       1& -q_1(x)
\\   0 & 1
  \end{matrix}\right),$$
and hence $\nu_{21}^{q^2}=\nu_{21}$, $\nu_{22}^{q^2}=\nu_{22}$ and
$\nu_{11}^{q^2}-\nu_{11}+q_1(x).\nu_{21}=0$,
 $\nu_{12}^{q^2}-\nu_{12}+q_1(x).\nu_{22}=0$. We can set
$\nu_{21}=\nu_{22}=1$ and $\nu_{11}-\nu_{12}= 1$.

In this case,
$$( \tilde\beta)^{-1}(Z^{\infty +0}_2) =\Delta$$
($\Delta$ denotes  the diagonal correspondence) because there exists a
morphism of level structures between two
vector bundles with the same degree and with $\infty   D$-level structures.
This morphism must be an isomorphism
(Proposition \ref{mor}).

As before, by   Theorem \ref{comp} and by the explicity calculations in
Remark \ref{matrix}  and Lemma
\ref{A} $  \tilde\beta $ has associated the $2\times
2$-matrix:
$$\left(\begin{matrix}  1  \otimes    \alpha^0_{1,0}& 1  \otimes
(\alpha^0_{1,0})^q
\\  1  \otimes  \alpha^0_{2,0}   &1  \otimes  (\alpha^0_{2,0})^q
  \end{matrix}\right)^{-1} .\left(\begin{matrix}      \alpha^0_{1,0}\otimes
1& (\alpha^0_{1,0})^q \otimes 1
\\   \alpha^0_{2,0}\otimes 1 & (\alpha^0_{2,0})^q\otimes 1
  \end{matrix}\right)
 \times $$
$$\left(\begin{matrix}     1  \otimes  \nu_{11} & 1  \otimes \nu_{12}
\\   1\otimes 1 & 1\otimes 1
  \end{matrix}\right)^{-1}.\left(\begin{matrix}       \nu_{11}\otimes 1&
\nu_{12}\otimes 1
\\   1\otimes 1 & 1\otimes 1
  \end{matrix}\right) $$

If we denote by
$$A^0_{0,1}:=\left(\begin{matrix}      \alpha^0_{1,0}\otimes 1&
(\alpha^0_{1,0})^q \otimes 1
\\   \alpha^0_{2,0}\otimes 1 & (\alpha^0_{2,0})^q\otimes 1
  \end{matrix}\right)\text{, }A^1_{0,1}:=\left(\begin{matrix}  1  \otimes
\alpha^0_{1,0}& 1  \otimes (\alpha^0_{1,0})^q
\\  1  \otimes  \alpha^0_{2,0}   &1  \otimes  (\alpha^0_{2,0})^q
  \end{matrix}\right)$$
and
$$A^0_\infty :=\left(\begin{matrix}       \nu_{11}\otimes 1& \nu_{12}\otimes 1
\\   1\otimes 1 & 1\otimes 1
  \end{matrix}\right)\text{, }A^1_\infty :=\left(\begin{matrix}     1
\otimes  \nu_{11} & 1  \otimes \nu_{12}
\\   1\otimes 1 & 1\otimes 1
  \end{matrix}\right) $$
then by Theorem \ref{comp}:
$$( \tilde \beta)^{-1}(Z^{\infty +0}_2)=( A^1_\infty (A^1_{0,1})^{-1}.
A^0_{0,1}.(A^0_\infty)^{-1}-Id_2=0)$$

\bigskip

{\bf Example. 4} Now,     $n=2$,
$p(t)=t(t-1) $. Thus, $D=0+1$ and as in the above example we follow  the
notation of   Remark \ref{matrix}:
$$\tilde K_{2}^{\infty +0 }={\Bbb F}_q(x)[ \alpha^0_{1,0}, \alpha^0_{2,0},
\beta^0_{1,0},
\beta^0_{2,0},\nu_{11},\nu_{21},\nu_{12},\nu_{22}],$$
where $\alpha^0_{i,0}$ is as in the above example and
$$(\beta^0_{i,0})^{q^2} +q_1(x).(\beta^0_{i,0})^q + \beta^0_{i,0}.(x-1)=0 $$
for $i=1,2$ and $(\nu_{ij})_{ij}$ as before.

$ A^0_\infty$, $A^1_\infty$,  $A^0_{0,1}$ and
$A^1_{0,1} $ follow  the same notation as in the last example
and $B^0_{0,1}$,
and $B^1_{0,1}$   denote
$$\left(\begin{matrix}      \beta^0_{1,0}\otimes 1& (\beta^0_{1,0})^q \otimes 1
\\   \beta^0_{2,0}\otimes 1 & (\beta^0_{2,0})^q\otimes 1
  \end{matrix}\right)\text{, }\left(\begin{matrix}  1  \otimes
\beta^0_{1,0}& 1  \otimes (\beta^0_{1,0})^q
\\  1  \otimes  \beta^0_{2,0}   &1  \otimes  (\beta^0_{2,0})^q
  \end{matrix}\right)$$
respectively.

Again by   Theorem \ref{comp} and by the explicit  calculations in   Remark
\ref{matrix}  and Lemma
\ref{A}, $\tilde \beta$ has
associated the $2\times 2$-matrix:
$$ (A^1_\infty)^{-1}. A^0_{0,1} \times (B^1_\infty)^{-1}. B^0_{0,1} \times
(A^1_\infty)^{-1}.
A^0_\infty$$
and   $( \tilde \beta)^{-1}(Z^{\infty +0+1}_2)$ is
 $$( A^1_\infty (A^1_{0,1})^{-1}.
A^0_{0,1}.(A^0_\infty)^{-1}-
A^1_\infty (B^1_{0,1})^{-1}.  B^0_{0,1}.(A^0_\infty)^{-1}-Id_2=0).$$
The coefficients $+$ and $-$ are deduced bearing in mind the
Chinese remainder theorem.
We shall now attempt to calculate   the irreducible components of $(
\tilde\beta)^{-1}(Z^{\infty
+0+1}_2)$. From Lemma \ref{beta}, we know that these components are
graphs $\Gamma_{g.F^i}$, where
$i\leq 2$ and $g\in Gl_2({\Bbb F}_q)_\infty\times Gl_2({\Bbb F}_q)_0\times
Gl_2({\Bbb F}_q)_1$.

\bigskip
We know from  Lemma \ref{irreducible} that the graphs $\Gamma_{F^2}$ and
$\Gamma_{h_{d+t}}$ are contained in
  $( \tilde\beta)^{-1}(Z^{\infty +0+1}_2)$.

We  try to prove that
$$( \tilde {\beta}_2^{\infty +0+1})^{-1}(\Theta^{\infty +0+1}_2)=\underset
{d\neq 0}\bigcup
\Gamma_{h_{d+t}}\cup \Gamma_{F^2}.$$
Let $g=g_\infty\times g_{0+1} $ be $B\times C\times D\in Gl_2({\Bbb
F}_q)_\infty\times Gl_2({\Bbb F}_q)_0\times
Gl_2({\Bbb F}_q)_1$.   If $\Gamma_g\subset ( \tilde\beta)^{-1}(Z^{\infty+0+1}_2)$,
we have
$$(A_{0,1})^{-1}.C.A_{0,1}-(B_{0,1})^{-1}.D.B_{0,1}=(A_\infty)^{-1}.B.A_\infty$$
where
$$B_{0,1}:=\left(\begin{matrix}      \beta^0_{1,0} & (\beta^0_{1,0})^q
\\   \beta^0_{2,0}  & (\beta^0_{2,0})^q
  \end{matrix}\right)
  A_{0,1}:=\left(\begin{matrix}      \alpha^0_{1,0} & (\alpha^0_{1,0})^q
\\   \alpha^0_{2,0}  & (\alpha^0_{2,0})^q
  \end{matrix}\right) $$
and
$$A_\infty :=\left(\begin{matrix}       \nu_{11} & \nu_{12}
\\   1  & 1
  \end{matrix}\right).$$
By acting on this equality by $Id\times \Omega\times Id \in Gl_2({\Bbb
F}_q)_\infty\times Gl_2({\Bbb F}_q)_0\times
Gl_2({\Bbb F}_q)_1$, we obtain the expression
$$(A_{0,1})^{-1}.\Omega^{-1}.C.\Omega.A_{0,1}-(B_{0,1})^{-1}.D.B_{0,1}=$$
$$=(A_\infty)^{-1}.B.A_\infty.$$
Thus, for each $\Omega\in Gl_2({\Bbb F}_q)$ we have that
$\Omega^{-1}.C.\Omega=C$, ans hence $C$ is a central matrix. In
an analogous way, we can check that $B$ and $D$ must be central matrices.
As ${\Bbb F}_q^*$ acts on $\tilde K_{2}^{\infty
+ 0+1}$ by the identity, we can assume that $B=Id$. Then, $C-D=Id$, and
hence $g=h_{d+t}$.

In a similar way, one can check that there does not exist $g$ with
$$\Gamma_{g.F}\subset  (\tilde \beta)^{-1}(Z^{\infty +0+1}_2).$$
We thus obtain what were looking for:
$$(\tilde \beta)^{-1}(Z^{\infty +0+1}_2)=\underset {d\neq 0}\bigcup
\Gamma_{h_{d+t}}\cup \Gamma_{F^2}.$$

\section{Case of a general curve}
In this  section we shall study zeta subschemes for a general curve $C$. We
shall translate the
results  of the above sections to this case. Now, instead of $\infty$ we
shall consider  a rational point $p\in
Spec(A)$  such that $p\notin Supp(D)$. In this case, we shall state in a
more general way the result of  Lemma
\ref{beta},
 and the relation between   these zeta subschemes and the general theta divisor.

Let us recall some notation and results:
${\cal M}_{C}(n,h)$ denotes the moduli stack, {\cite {LM}},  of
 vector bundles of rank $n$ and degree $h$. In the same way,
${\cal M}_{C}(n,h,p+ D)$ denotes the moduli stack of pairs, $(M,f_{p+
D})$, of
 vector bundles of rank $n$ and degree $h$ with a ${p+   D}$-level structure.
 There exists a natural epimorphism of forgetting the level structures:
$$L: {\cal M}_{C}(n,h,p+ D)\to{\cal M}_{C}(n,h)$$
Let us denote  by $\Theta_n  \subset {\cal M}_{C}(n,n(g-1))$ the general
theta divisor \cite{DN}:
$\Theta_n$ is the closed substack of ${\cal M}_{C}(n,n(g-1))$ of vector
bundles $M$ such that
$H^0(C,M)=0$. We set $\Theta_n^{\infty   D }:=L^{-1}(\Theta_n) $.

By ${\cal E}_n^{p+  D}$ we denote the moduli (scheme for $deg(D)>>0$) of
rank $n$-$A$-elliptic sheaves, $(E_j,i_j,\tau)$,
 with   $p+D$-level structures.
As in the preceding sections, we have a natural morphism
$$\Phi:  {\cal E}_n^{p+   D}\to {\cal M}_{C}(n,ng,p+ D),$$
given by $\Phi(E_j,i_j,\tau, f_{p+D})=(E_0, f^0_{p+D})$.

In this section,  $\pi$ denotes the natural morphism
$\o_C(g+d-1)\to \o_C/\o_C(-p-D)$. As in the above sections
$$\o_{p D}:=\o_C/\o_C(-p-D)\text{ and }O_{p
D}:=H^0(C,\o_C/\o_C(-p-D)).$$

\begin{rem}
We consider  the open substack
$${\cal U}_d\subseteq {\cal M}_{C}(n,h,p+ D)\times {\cal M}_{C}(n,h,p+ D)$$
of pairs $[(\bar M,\bar f_{p+  D}),(M,f_{p+   D})]$ satisfying $H^1(C,{\bar
M}^\vee \underset{\o_C } \otimes
M(g+d-1))=0$. In this case, if $[(\bar {\Bbb M},\bar f_{p+  D}),({\Bbb
M},f_{p+   D})]$ is a  universal object for
${\cal U}_d$ we have
$$h^0({\bar {\Bbb M}}^\vee \underset{\o_{C \times {\cal U}_d}} \otimes
{\Bbb M}(g+d-1))=n^2d.$$

Let
$\{s_1,\cdots,s_{n^2d}\}$ be a basis for the space global sections of
$${\bar {\Bbb M}}^\vee \underset{\o_{C \times {\cal
U}_d} } \otimes {\Bbb M}(g+d-1),$$
and we set
$$o_i=H^0(\bar f_{p+  D}^\vee\otimes f_{p+ D}\otimes\pi)(s_i).$$

\end{rem}

In this last Remark and in   the following Lemmas we assume, for easy
notation, that the global sections of
${\bar {\Bbb M}}^\vee \underset{\o_{C \times {\cal U}_d}} \otimes {\Bbb
M}(g+d-1)$ are a free module  of rank
$n^2d$. By Grauert's Theorem, it occurs locally, so instead of ${\cal U}_d$
we should write ${{\cal U}_d}^\alpha$,  where
$ {\cal U}_d=\underset \alpha \cup {\cal U}_d^\alpha$. The localization by
${\cal U}_d^\alpha$ of the space of
global sections of
$${\bar {\Bbb M}}^\vee \underset{\o_C \otimes {\o_{\cal U}}_d } \otimes
{\Bbb M}(g+d-1)$$
 is a
$\o_{{\cal U}_d^\alpha}$-free module.

 We shall now define an analogous subscheme to $Z^{\infty  D}_n$.

\begin{defn}\label{zeta2} Let $S$ be a scheme. We define the $p+
D$-generalized zeta substack, $Z^{p+   D}_n$,
as the closed substack of   ${\cal U}_d$ defined by:
$[(\bar M,\bar f_{p+  D}),(M,f_{p+   D})]\in Hom_{stacks}(S,Z^{p+   D}_n) $
if  over the open subscheme, $\bar U$, of $S$
formed by the $s\in S$ such that
$$H^0(C\otimes k(s), {\bar  M}_s^\vee \underset{\o_{C \times k(s)}} \otimes
M_s((g+d-1)\infty-p-D))=0$$
  there exists a morphism of modules
$T:\bar M_{\vert \bar U}\to M(g+d-1)_{\vert \bar U} $ ($d=deg(D)$), where
the diagram :

 $$ \xymatrix {  \bar M_{\vert \bar U}   \ar[dr]_{\bar f_{p+   D}}
\ar[r]^(.3){T} &
  M(g+d-1)_{\vert \bar U}\ar[d]^{f_{p+  D} \otimes \pi}
\\ &(\o_{p D})^n \otimes \o_{\bar U} } $$
is commutative.

\end{defn}

\begin{rem} From the last Lemma, if
$$[(\bar M,\bar f_{p+  D}),(M,f_{p+   D})]\in Hom_{stacks}(S,Z^{p+   D}_n) $$
then there exists an open subscheme $\bar U\subset S$ such that over this
open subscheme there exists a non-trivial
morphism
$T:\bar M
\to M(g+d-1)
$.   And if $s\notin  \bar U$,
there exists a non-trivial   morphism
$T_s:\bar M_s\to M_s((g+d-1)\infty-p-D)) \subset M(g+d-1)_s $.

Hence $Z^{p+   D}_n=V^{p+   D}_n\cup (V^{p+   D}_n)^c$, with $V^{p+   D}_n$
an open substack of $Z^{p+   D}_n$.

\end{rem}

\begin{rem}
A pair $[(\bar M,\bar f_{p+  D}),(M,f_{p+   D})]\in
Hom_{stacks}(Spec(R),Z^{p+   D}_n)
$  gives a morphism
$$  \bar f_{p+  D}^\vee\otimes f_{p+
D}\otimes \pi:  {\bar M}^\vee \underset{\o_C\otimes R } \otimes M(g+d-1)
 \longrightarrow (\o_{pD})^{n^2}\otimes R. $$
Therefore, by taking global sections we obtain a morphism
$$H^0(\bar f_{p+  D}^\vee\otimes f_{p+ D}\otimes \pi): H^0(\o_C\otimes
R,{\bar M}^\vee
\underset{\o_C\otimes R } \otimes M(g+d-1)) \to (O_{pD})^{n^2}\otimes R.$$

We shall denote by $\Xi$ the natural morphism
$$\Xi:\o_C\otimes R \to {\bar M}^\vee \underset{\o_C\otimes R } \otimes \bar M
\overset{\bar f_{p+  D}^\vee\otimes \bar f_{p+
D} }\longrightarrow (\o_{pD})^{n^2}\otimes R $$
where the first morphism in the composition is the diagonal morphism.
Therefore, by taking global sections we obtain a
morphism
$$H^0(\Xi):R\to  (O_{pD})^{n^2}\otimes R$$
given by the diagonal matrix with entries over $O_{pD}$.

Let us now consider  the morphism
$$\Sigma:=H^0(\Xi)+H^0(\bar f_{p+  D}^\vee\otimes f_{p+ D}\otimes \pi)$$
from
$R\oplus H^0(\o_C\otimes R,{\bar M}^\vee \underset{\o_C\otimes R } \otimes
M(g+d-1))$ to $(O_{pD})^{n^2}\otimes R$.
\end{rem}

\begin{lem}\label{Sigma} A pair $[(\bar M,\bar f_{p+  D}),(M,f_{p+
D})]\in Hom_{stacks}(Spec(R),{\cal U}_d)$ belongs to
$$Hom_{stacks}(Spec(R),Z^{p+   D}_n)$$
if and only if $\Sigma$ is not injective.
\end{lem}
\begin{proof} Since   if $s\in Spec(R)$ satisfies
$$H^0(C\otimes k(s), {\bar  M}_s^\vee \underset{\o_{C \times k(s)}} \otimes
M_s((g+d-1)\infty-p-D))\neq 0$$
 we have  $Ker(H^0(\bar f_{p+  D}^\vee\otimes f_{p+ D}\otimes \pi)_s)\neq
0$. Thus, $\Sigma_s$ is not injective and we can
assume    that for all   $s\in Spec(R)$ we have:
$$H^0(C\otimes k(s), {\bar  M}_s^\vee \underset{\o_{C \times k(s)}} \otimes
M_s((g+d-1)\infty-p-D))= 0.$$

  If there exists a morphism
$T$ between the  level structures  considered
$$ \xymatrix {  \bar M   \ar[dr]^{\bar f_{p+   D}}    \ar[r]^(.3){  T} &
  M(g+d-1)\ar[d]^{f_{p+  D} \otimes \pi}
\\ &(\o_{p D})^n \otimes R } $$
by tensoring by $\bar M^\vee$ and considering $\bar f_{p+  D}$, we obtain a
section
 $$w\in H^0(\o_C\otimes R,{\bar M}^\vee \underset{\o_C\otimes R } \otimes
M(g+d-1))$$
such that
$$H^0(\Xi)(1)=H^0(\bar f_{p+  D}^\vee\otimes f_{p+ D}\otimes \pi)(w).$$
Thus, in this case $1\oplus-w\in Ker(\Sigma)$.

Conversely, if $1\oplus-s\in Ker(\Sigma)$ we have a commutative diagram
$$ \xymatrix {\o_C\otimes R   \ar[dr]_\Xi     \ar[r]  &
  {\bar M}^\vee \underset{\o_C\otimes R } \otimes M(g+d-1)\ar[d]^{\bar
f_{p+  D}^\vee\otimes f_{p+ D}\otimes
\pi}
\\ &(\o_{p D})^{n^2} \otimes R } .$$
Therefore, we conclude by tensoring by $\bar M$ and bearing in mind $\bar
f_{p+  D}$.

\end{proof}

\begin{rem}  If $Z^{p+   D}_n=V^{p+   D}_n\cup (V^{p+   D}_n)^c$, we deduce
that   that:
$$(V^{p+   D}_n)^c =\{[(\bar M,\bar f_{p+  D}),(M,f_{p+   D})]\in Z^{p+
D}_n\text{, with }
  Ker (H^0(\bar f_{p+  D}^\vee\otimes f_{p+ D}\otimes \pi) \neq 0 \}$$
and
 $[(\bar M,\bar f_{p+  D}),(M,f_{p+   D})] \in V^{p+   D}_n$ if
$$H^0(\Xi)(1)\in H^0(\bar f_{p+  D}^\vee\otimes f_{p+ D}\otimes \pi).$$
\end{rem}

\begin{rem}
From $H^0(\Xi)(1)$ and $\{o_1,\cdots,o_{n^2d}\}$, we obtain a $
(n^2d+1)(columns)\times (n^2(d+1))(rows)$-matrix
$(H^0(\Xi)(1), o_1,\cdots,o_{n^2d})$ with  entries over $H^0({\cal
U}_d,{\o_{\cal U}}_d)$.  Recall that $ \{H^0(\Xi)(1),
o_1,\cdots,o_{n^2d}\}\subset ( O_{pD})^n\otimes H^0({\cal U}_d,{\o_{\cal
U}}_d)$ and $dim_{{\Bbb F}_q}O_{pD}=d+1$.
\end{rem}

Let $\bar {\cal U}_d$ be the open substack of ${\cal U}_d$ such that
$[(\bar M,\bar f_{p+  D}),(M,f_{p+   D})]
\in \bar {\cal U}_d$ if and only if $h^0({\bar {\Bbb M}}^\vee
\underset{\o_C } \otimes
{\Bbb M}(g+d-1)(-p-D))=0$. In this case,  $rank (  o_1,\cdots,o_{n^2d})=n^2d $.

\begin{lem} $ Z^{p+ D}_n$ is a closed substack of ${\cal U}_d$ and   is
formed by the pairs
$$[(\bar M,\bar f_{p+  D}),(M,f_{p+   D})]\in {\cal U}_d,$$
such that
$$rank ((O_{p+D})^{n^2}/<H^0(\Xi)(1), o_1,\cdots,o_{n^2d} >)\geq n^2 .$$

Moreover, $V^{p+   D}_n=\bar {\cal U}_d\cap Z^{p+   D}_n$
is locally the zero locus of $n^2$-regular functions.
\end{lem}

\begin{proof}The first assertion is easily deduced from Lemma \ref{Sigma}.

For the second assertion,
we denote by
$$N^{(H^0(\Xi)(1), o_1,\cdots,o_{n^2d})}_{i_1,\cdots,i_{n^2d+1}}$$
($1\leq i_1<\cdots< i_{n^2d+1}\leq n^2(d+1)$)
the $n^2d+1$-minor of the
matrix
$$(H^0(\Xi)(1), o_1,\cdots,o_{n^2d}),$$
 where the $k^{th}$-row is the ${i_k}^{th}$-row of $(H^0(\Xi)(1),
o_1,\cdots,o_{n^2d})$.

Analogously, we set
$$N^{(  o_1,\cdots,o_{n^2d})}_{i_1,\cdots,i_{n^2d }}$$
for the
$ (n^2d )(columns)\times(n^2(d+1))(rows)$-matrix $(  o_1,\cdots,o_{n^2d})$.
Let us take the open covering of $\bar {\cal
U}_d$
$$\bar {\cal U}_d=\underset{1\leq i_1<\cdots<
i_{n^2d  }\leq n^2(d+1)} \bigcup \bar {{\cal U}_d}^{ i_1,\cdots, i_{n^2d }}$$
 $\bar{{\cal U}_d}^{i_1,\cdots, i_{n^2d }}$ being   the open substack of
$\bar {\cal U}_d$ of pairs $[(\bar M,\bar f_{p+
D}),(M,f_{p+   D})]\in {\cal U}_d$, such that
$$rank( N^{(  o_1,\cdots,o_{n^2d})}_{i_1,\cdots,i_{n^2d }})=rank (
o_1,\cdots,o_{n^2d})=n^2d.$$
Thus, we have that $\bar {{\cal U}_d}^{i_1,\cdots, i_{n^2d }}\cap Z^{p+
D}_n$ is the zero locus of the functions
$$\{det(N^{(H^0(\Xi)(1), o_1,\cdots,o_{n^2d})}_{i_1,\cdots,i_{n^2d },
j_{n^2d+1} }),\cdots,
det(N^{(H^0(\Xi)(1), o_1,\cdots,o_{n^2d})}_{i_1,\cdots,i_{n^2d },
j_{n^2d+n^2} })\},$$
with $\{1,2,\cdots,n^2(d+1)\} \setminus \{i_1,\cdots,i_{n^2d
}\}=\{j_{n^2d+1},\cdots,j_{n^2d+n^2}\}$.

\end{proof}
If
$$[(\bar M,\bar f_{p+  D}),(M,f_{p+   D})]
\in  Hom_{stacks}(Spec(R),\bar {\cal U}_d),$$  we have
$$H^0(C\otimes R,{\bar  M}^\vee \underset{\o_C \otimes R} \otimes
 M(g+d-1)(-p-D))=0$$ and
$$H^1(C\otimes R,{\bar  M}^\vee \underset{\o_C \otimes R} \otimes
 M(g+d-1) )=0.$$
In this way,
$$H^0(\bar f_{p+  D}^\vee\otimes f_{p+ D}\otimes
\pi)(H^0(C\otimes R,{\bar  M}^\vee \underset{\o_C } \otimes
 M(g+d-1) ))$$
 is a rank $n^2d$-sub-vector bundle of $(O_{pD})^{n^2}\otimes R $. Thus,
$$\overset{n^2d }\bigwedge (H^0(\bar f_{p+  D}^\vee\otimes f_{p+ D}\otimes
\pi)(H^0( {\bar  M}^\vee \underset{\o_C \otimes R} \otimes
 M(g+d-1) )) \in {\Bbb P}(\overset{n^2d }\bigwedge (O_{pD})^{n^2}\otimes R) $$
and we obtain a morphism
$$\bar {\cal U}_d \to {\Bbb P}(\overset{n^2d }\bigwedge (O_{pD})^{n^2}).$$

It is not hard to prove that
$$det(N^{(H^0(\Xi)(1), o_1,\cdots,o_{n^2d})}_{i_1,\cdots,i_{n^2d },
j_{n^2d+k} })$$
can be considered as  the restriction to $\bar {\cal U}_d$ of a certain
hyperplane  section
  of
${\Bbb P}(\overset{n^2d }\bigwedge (O_{pD})^{n^2})$. In this way, ${{\cal
U}_d}^{i_1,\cdots, i_{n^2d }}\cap Z^{p+   D}_n$
is obtained as the intersection  of $n^2$ hyperplane sections  of ${\Bbb
P}(\overset{n^2d }\bigwedge (O_{pD})^{n^2})$.

\bigskip
For $n=1$, we have a direct relation of $Z^{p+   D}_1$ with a theta
divisor. We take $d>2g-2$.

\begin{defn} We define the $p+  D$-generalized theta divisor, $\Theta_{p+
D } $,
as the closed subscheme of   $J^{p+  D}_{C,0} $ defined by:
$(L,f_{p+   D})\in Hom_{schemes}(Spec(K),\Theta_{p+  D} ) $ ($K$ a field)
if and only if there exists a morphism of
modules
$T:\o_{C} \otimes R\to L(g+d-1) $ ($d=deg(D)$), where the diagram:

$$ \xymatrix {  \o_{C} \otimes K\ar[d] \ar[r]^{T} & L(g+d-1)\ar[dl]^{f_{p+
D}\otimes
\pi } \\
\o_{p+    D}\otimes K & }$$
is commutative. The vertical arrow is the natural epimorphism.
\end{defn}

Let us consider    $J^{p+   D}_{C,g}$, where $J^{p+   D}_{C,g} $ is the $p+
D$-generalized Rosenlich's Jacobian
(c.f:\cite{Se}).
\begin{rem} By regarding the morphism
$$m_{p+ D}: J^{p+   D}_{C,g}  \times
 J^{\infty   D}_{C,g}  \to  J^{p+  D}_{C,0}  $$
 defined by
$$m_{p+ D}[( \bar L,  \bar f_{p+   D}),(  L,  f_{p+   D})]=( \bar
L^\vee,\bar f_{p+   D}^\vee)\otimes( L ,
 f_{p+   D}),$$
we have that $m_{\infty D}^{-1}(\Theta^{\infty   D}_1 )=Z^{\infty   D}_1$
for $d>2g-2$, because of the definition
of $Z^{\infty   D}_1$ and
${\cal U}_d=J^{p+   D}_{C,g}
\times
 J^{\infty   D}_{C,g}$.

Moreover, from the last Lemma we can say  that $\bar {{\cal
U}_d}^{i_1,\cdots, i_d }\cap Z^{p+   D}_1$ is the zero locus
of the function
$$det(N^{(1, o_1,\cdots,o_d)}_{i_1,\cdots,i_d ,  j_{d+1} })$$
over $\bar {{\cal U}_d}^{i_1,\cdots, i_d }$ (for $n=1$  $H^0(\Xi)(1)=1\in
O_{pD}$). For different
$\{i_1,\cdots,i_d
\}$, these determinants are the determinant
$det(1,
o_1,\cdots,o_d)$
  up to a sign. Then,
  $$\bar {{\cal U}_d}\cap Z^{p+   D}_1=\bar {{\cal U}_d}\cap (det(1,
o_1,\cdots,o_d))_0.$$
Recall that
$$\bar {\cal U}_d=\underset{1\leq i_1<\cdots<
i_{n^2d  }\leq d+1} \bigcup \bar {{\cal U}_d}^{ i_1,\cdots, i_2d }.$$
On the other hand, as $(\bar {\cal U}_d)^c$ is formed by the pairs where
$\{o_1,\cdots,o_d\}$ are linearly dependent,
then
$$(\bar {\cal U}_d)^c\subseteq (det(1, o_1,\cdots,o_d))_0,$$
  and we conclude that $  Z^{p+
D}_1=(det(1, o_1,\cdots,o_d))_0$.
\end{rem}

Recall that actually this result is settled locally over $J^{p+
D}_{C,g}\times J^{p+   D}_{C,g}$.   To be precise, one
must take a universal object
$[(\bar {\Bbb L}, \bar f),(  {\Bbb L},  f)]$
and an open covering for $J^{p+   D}_{C,g}\times J^{p+   D}_{C,g}$  such that
$$p_*(\bar {\Bbb L}^\vee \underset{\o_C\otimes \o_{J}}\otimes  {\Bbb L})$$
is trivialized by this covering. $p$ denotes the natural projection
$C\times J^{p+   D}_{C,g}\to J^{p+   D}_{C,g}$.

\subsection{Zeta subschemes for rank-$n$-Shtukas}

In this section we shall study     zeta subschemes in the case of
rank-$n$-Shtukas and $A$-Drinfeld modules defined over
an arbitrary curve $C$. Professor G.W. Anderson has suggested that I should
study this topic this topic.

First, we shall state some previous definitions.

\begin{defn} A rank-$n$ shtuka over $C\times Spec(R)$ is a diagram  of
rank-$n$-vector bundles over $C\times Spec(R)$
$${\cal L}_r\overset{j} \hookrightarrow {\cal L}_{r+1} \overset{t}
\hookleftarrow F^{\#}{\cal L}_r,$$
such that ${\cal L}_{r+1} /j({\cal L}_r)$ and ${\cal L}_{r+1}/t(F^{\#}{\cal
L}_r)$   are line bundles as $R$-modules and
they are supported by two morphisms
$\Gamma_\alpha$ and
$\Gamma_\gamma$  $Spec(R)\to C$; the pole and zero, respectively. We denote the
line bundles associated with $\Gamma_\alpha$ and
$\Gamma_\gamma$ by $\o_{C\otimes R}(\alpha)$ and $\o_{C\otimes R}(\gamma)$
respectively. We assume
 $j({\cal L}_r)+t(F^{\#}{\cal L}_r)={\cal L}_{r+1}$. We set $deg({\cal
L}_r)=ng+r$.
 \end{defn}

\begin{defn} Similar to the above sections, to give a $p+D$-level structure
over
$${\cal L}_r\overset{j} \hookrightarrow {\cal L}_{r+1} \overset{t}
\hookleftarrow F^{\#}{\cal L}_r$$
is to give level structures $({\cal L}_r,f^r_{p+D})$ and $({\cal
L}_{r+1},f^{r+1}_{p+D})$ with
commutative diagrams
$$ \xymatrix {  {\cal L}_r   \ar[dr]^{f^r_{p+D}}    \ar[r]^{j} &
 {\cal L}_{r+1}\ar[d]^{f^{r+1}_{p+D}}
\\ &(\o_{p+ D})^n \otimes R } $$
 and
$$ \xymatrix {   F^{\#}{\cal L}_r   \ar[dr]^{ F^{\#} f^r_{p+D}}    \ar[r]^{t} &
 {\cal L}_{r+1}\ar[d]^{f^{r+1}_{p+D}}
\\ &(\o_{p+ D})^n \otimes R }, $$
$p\notin \{\alpha, \gamma, supp(D)\}$.
\end{defn}

\begin{defn} An isogeny between two rank-$n$ Shtukas
$${\cal L}_r\overset{j} \hookrightarrow {\cal L}_{r+1} \overset{t}
\hookleftarrow F^{\#}{\cal L}_r$$
and
$${\cal N}_h\overset{j'} \hookrightarrow {\cal N}_{h+1} \overset{t'}
\hookleftarrow F^{\#}{\cal N}_h$$
is defined by two morphisms of $\o_{C\otimes R}$-modules $T: {\cal L}_r\to
{\cal N}_h$ and $H : {\cal L}_{r+1}\to {\cal
N}_{h+1}$ such that  $j'.T=H.j$ and  $t'.F^{\#} T=F^{\#} H.t$.

 Analogously, one can define a  $p+D$-isogeny between shtukas with
$p+D$-level structures,  except that one must take $T$
and
$H$ compatible  with the level structures.
\end{defn}

We now consider  the stack, $Cht^{n,ng+r}$, of shtukas of rank $n$  and
$deg({\cal L}_r)=ng+r$.
We follow  the same notation for $Cht_{p+D}^{n,ng+r}$   but considering
$p+D$-level structures. Let
$P^{n,ng+r }$ be the closed substack of $ Cht^{n,ng+r}\times Cht^{n,ng+r}$
of shtukas with the same pole. Analogously, we
define
$P_{p+D}^{n,ng+r }$ .

\begin{propo}\label{isosh} Let $[({\cal L}_r,{\cal L}_{r+1},j,t),({\cal
N}_h,{\cal N}_{h+1},j',t')]$ be an  element of $
P^{n,ng+r }$ with pole
$\alpha$ and defined over a
 ${\Bbb F}_q$-reduced algebra $R$.
Then, there exists an isogeny between these shtukas if and only if there
exists a morphism of modules
$T:{\cal L}_r\to {\cal N}_h,$ together with the commutative diagram
$$ \xymatrix {  {\cal L}_r(\alpha)  \ar[r]^{T(\alpha)}    &
 {\cal N}_h (\alpha)
\\ {\cal L}_{r+1}\ar@{^{(}->}[u]^\kappa   &{\cal
N}_{h+1}\ar@{^{(}->}[u]^{\kappa'}
\\ F^{\#}{\cal L}_r\ar@{^{(}->}[u]^{t}\ar[r]^{F^{\#} T }   &F^{\#}{\cal
N}_h\ar@{^{(}->}[u]^{t'}    } $$
  $T(\alpha)$ denotes the morphism induced by $T$ from ${\cal L}_r(\alpha)$ to
 ${\cal N}_h (\alpha)$. The injective morphisms $ \kappa:{\cal
L}_{r+1}\hookrightarrow {\cal L}_r(\alpha)$ and
$ \kappa':{\cal N}_{h+1}\hookrightarrow  {\cal N}_h(\alpha)$ are  morphisms
such that the composition with $j$ and $j'$ are
the morphisms
${\cal L}_r\hookrightarrow {\cal L}_r(\alpha)$ and ${\cal
N}_h\hookrightarrow  {\cal N}_h(\alpha)$ induced by the natural inclusion
$\o_{C}\otimes R\subset \o_{C}\otimes R(\alpha)$.
Recall that the cokernels of
${\cal L}_r\hookrightarrow {\cal L}_{r+1}$ and ${\cal
N}_h\hookrightarrow  {\cal N}_{h+1}$   are supported over $\alpha$.
\end{propo}
\begin{proof}
This can be deduced because
$$ \xymatrix {  {\cal L}_r(\alpha)  \ar[r]^{T(\alpha)}    &
 {\cal N}_h (\alpha)
\\ {\cal L}_{r+1}\ar@{^{(}->}[u]^{\kappa }    &{\cal
N}_{h+1}\ar@{^{(}->}[u]^{\kappa'}
\\ {\cal L}_r\ar@{^{(}->}[u]^j\ar[r]^{ T }   &{\cal
N}_h\ar@{^{(}->}[u]^{j'}    } $$
is commutative, and
$$j({\cal L}_r)+t( F^{\#}{\cal L}_r)= {\cal L}_{r+1}  \text{, }j'({\cal
N}_h)+t'( F^{\#}{\cal N}_h)= {\cal N}_{h+1}.$$
\end{proof}

Let consider us the natural morphism
$\bar \Phi:Cht_{p+D}^{n,ng+r} \to {\cal M}_{C}(n,ng+r,p+ D)$
given by $\bar \Phi(({\cal L}_r,f^r_{p+D}),({\cal
L}_{r+1},f^{r+1}_{p+D}),j,t))=({\cal L}_r,f^r_{p+D})$.
Hence, we deduce a morphism
$$\bar \beta: P_{p+D}^{n,ng+r }\to Cht_{p+D}^{n,ng+r}\times
Cht_{p+D}^{n,ng+r} \overset {\bar \Phi\times \bar
\Phi}\longrightarrow {\cal M}_{C}(n,ng+r,p+ D)\times {\cal M}_{C}(n,ng+r,p+
D).$$
For easy notation, we denote by $({\cal L}_r,{\cal L}_{r+1},f_{p+D})$ an
element  of $Cht_{p+D}^{n,ng+r}$. Because
$p\neq
\alpha$ in the following Lemma, we are be able to set, for shtukas, the  results of   Lemma
\ref {cuadrado}   in a more precise way . This Lemma attempts to be an
approximation for shtukas with level structures  of Lemma 3.3 \cite{An2}.

\begin{lem}\label{isogenia}  Let $[({\cal L}_0,{\cal L}_{1},f_{p+D}),
({\cal N}_0,{\cal N}_{1},g_{p+D})] \in \bar
\beta^{-1}(Z^{p+   D}_n)$   be defined over a field
$K$.  Then:

 1) There exists a  $p+D$-isogeny between the shtukas with $p+D$-level
structures $[({\cal L}_0,{\cal L}_{1},f_{p+D})$
and $({\cal N}_0,{\cal N}_{1},g_{p+D}\otimes\pi)]$.

or

2) There exists an isogeny between the shtuka  $ ({\cal L}_0,{\cal
L}_{1},j,t)$ and the twisted shtuka
 $({\cal N}_0   ,{\cal
N}_{1} ,j',t')\underset{\o_C}\otimes \o_C((g+d-1)\infty-p-D)$.

or

3) $(F^{\#}{\cal L}_0)^\vee\underset {\o_C\otimes K}\otimes {\cal N}_0
((g+d-1)\infty  -p-D+\alpha)\in \Theta_{n^2}$.
\end{lem}
\begin{proof} Let us consider   the diagram  (not necessarily commutative)

$$ \xymatrix { {\cal L}_0(\alpha)   \ar[rr]^{T(\alpha)} \ar[dr]      & &
 {\cal N}_0(g+d-1) (\alpha) \ar[dl]
\\ {\cal L}_{1}\ar[u]^{\kappa }&(\o_{p+ D})^n\otimes R &{\cal
N}_{1}(g+d-1)\ar[u]^{\kappa' }
\\ F^{\#}{\cal L}_0\ar[u]^{t }\ar[rr]^{F^{\#} T } \ar[ur] &  &F^\#{\cal
N}_0(g+d-1)\ar[u]^{t' }
\ar[ul]    }
$$
$T$ is a morphism given by $[({\cal L}_0,f^0_{p+D}), ({\cal
N}_0,g^0_{p+D})] \in Z^{p+   D}_n$. The diagonal
morphisms are the level structures and the vertical morphisms $\kappa$ and
$  \kappa'$ are induced by
$j$ and $j'  $ since these two  shtukas have the same pole, $\alpha$. From
Definition \ref{zeta2}, we have two
possibilities:  either $T:{\cal L}_0   \to
 {\cal N}_0(g+d-1) $ takes values in  ${\cal N}_0(g+d-1)(-p-D)$ or the
latter diagram  is commutative after
tensoring by $\o_{pD}$. In both cases, one obtains a morphism:
$$T(\alpha).\kappa.t-\kappa'.t'.F^{\#} T:F^{\#}{\cal L}_0 \to {\cal
N}_0(g+d-1) (-p-D+\alpha).$$
 When $T(\alpha).\kappa.t-\kappa'.t'.F^{\#} T=0$, we obtain
  cases 1) and 2) and when it is $\neq 0$ it tell us that $(F^{\#}{\cal
L}_0)^\vee\underset {\o_C\otimes R}\otimes {\cal
N}_0 ((g+d-1)\infty  -p-D+\alpha)\in \Theta_{n^2}$.
\end{proof}


\subsection{ Zeta subschemes and $A$-Drinfeld modules}

 In this subsection, by means of the natural immersion of elliptic sheaves
into shtukas (c.f.\cite{Dr3}), we shall
be able to  translate the results of the above section to the case of
Drinfeld modules. We use  the notation and
definitions of section 3.

\begin{propo} There exists an immersion of stacks
$$\Upsilon:{\cal E}_n ^{p+   D}\hookrightarrow Cht_{p+D}^{n,ng}$$
\end{propo}
\begin{proof} This is defined by
$$\Upsilon(E_j,i_j,\tau,f_{p+D})=((E_0,f^0_{p+D}),
(E_1,f^1_{p+D}),i_{0},\tau)$$
This immersion is given in \cite{Dr3} (1.3).
\end{proof}

\bigskip
We recall the definition of isogeny for $A$-Drinfeld modules
\begin{defn} An isogeny of degree $r$ between two Drinfeld modules $\phi$
and $\bar \phi$ is a polynomial of degree $  r$
 in $\sigma$, $q(\sigma)$, such that the endomorphism $q(\sigma):({\Bbb
G}_a)_R\to ({\Bbb G}_a)_R$ satisfies
$$q(\sigma)(\phi_a(\lambda))=\bar \phi_a (q(\sigma))(\lambda))$$
for all $a\in A$.
\end{defn}

\begin{defn}Two elliptic sheaves  $(\bar E_{j},\bar i_{j},\bar \tau)$ and
$(  E_{j},   i_{j},  \tau)$ are said to be isogenous by an   isogeny of
degree $\leq r$ if and only if there exist morphisms of modules $T_j:\bar
E_{j}\to   E_{j+r}$, for each $j$, compatible
with the diagrams that define the elliptic sheaves.

Analogously, one can define a  $p+D$-isogeny between elliptic sheaves with
$p+D$-level structures, except that one must
take "$T_j$"  compatible  with the level structures.
\end{defn}

There exists an equivalence  between the  notion of isogenous elliptic
sheaves and isogenous   Drinfeld modules:
\begin{propo} Two Drinfeld modules $\phi$ and $\bar \phi$   associated with
elliptic sheaves  $(\bar E_{j},\bar
i_{j},\bar \tau)$ and
$(  E_{j},   i_{j},  \tau)$, respectively,
 are isogenous by an isogeny of
degree
$\leq r$ if and only if there exists an isogeny of
degree
$\leq r$ between the elliptic sheaves
$  (\bar E_{j},\bar i_{j},\bar \tau)$ and
 $(  E_{j+r},   i_{j+r},  \tau)$.
\end{propo}

The following Proposition together with   Proposition \ref{isosh} tell us
that two elliptic sheaves are isogenous as
shtukas if and only if they are isogenous as
 elliptic sheaves. Recall that  the   elliptic sheaves considered have
$\infty $ as their pole.

\begin{propo} Two elliptic sheaves  $(\bar E_j,\bar i_j,\bar \tau)$ and
$( E_j,  i_j, \tau)$ defined over a ${\Bbb F}_q$- reduced algebra, $R$,
 are isogenous by an isogeny of
degree
$\leq r$ if and only if there exists a morphism of vector bundles
$$T: \bar E_0\to    E_{r}$$
with the diagram
$$ \xymatrix { \bar E_0  \ar[rr]^{T}     & &
   E_{r}
\\ F^\#\bar E_{-n}\ar[u]^{\bar \tau}\ar[rr]^{F^\#T(-1)}  & & F^\#
E_{r-n}\ar[u]^{  \tau}   } $$
commutative. Here, $T(-1)$ denotes the morphism induced by $T$ over
$$\bar E_{-n}\to   E_{r-n},$$
bearing in mind that $E_{k-n}:=E_k(-\infty)$.
\end{propo}

\begin{proof} Since $\bar E_{1-n}=\bar\tau(F^\#\bar E_{-n})+\bar
E_{-n}\subset \bar E_0$ and $  E_{r+1-n}=  \tau(
F^\# E_{r-n} )+  E_{r-n}\subset E_r$, by observing the latter diagram  we
have $T(\bar E_{1-n})\subset   E_{r+1-n}$.
Hence, we obtain the commutative diagram

$$ \xymatrix { \bar E_0  \ar[r]^{T}      &
   E_{r}
\\\bar E_{1-n}\ar@{^{(}->}[u]\ar[r]^{T_{1-n}}   &  E_{r+1-n}\ar@{^{(}->}[u]
\\ \bar E_{-n}\ar@{^{(}->}[u]\ar[r]^{T(-1)}   &\  E_{r-n}\ar@{^{(}->}[u]
}, $$

where $T_{1-n}$ denotes the restriction of $T$ to $\bar E_{1-n}$  and the
injective morphisms are the "$i$" and "$\bar i$"
inclusions given in the elliptic sheaves.

From the last commutative diagram  we also obtain a commutative diagram
$$ \xymatrix {  F^\#\bar E_{1-n} \ar[rr]^{F^\# T _{1-n}} &  & F^\# E_{r+1-n}
\\F^\#\bar E_{-n}\ar@{^{(}->}[u]\ar[rr]^{F^\#T(-1)}  & &
F^\#E_{r-n}\ar@{^{(}->}[u]    } $$
 and in the case of $R$  being a field,
$$ \xymatrix { \bar E_0      \ar[rr]^{T} & &
  E_{r}
\\ F^\#\bar E_{1-n}\ar[u]^{\bar\tau}\ar[rr]^{F^\#T_{1-n}} &  & F^\#
E_{r+1-n}\ar[u]^{  \tau}   }, $$

($\bar \tau$ and $\tau$ are morphisms induced by the definition of elliptic
sheaves) is also commutative because $(T.\bar
\tau-
\tau.F^\#T_{1-n})_{\vert
F^\#\bar E_{-n}}=0
$. Then,
$$F^\#\bar E_{-n}\subset Ker (T.\bar\tau- \tau.F^\#T_{1-n})$$
and hence rank($Ker
({T.\bar \tau- \tau.F^\#T_{1-n}}  ))=n$. Because
 $T.\bar\tau- \tau.F^\#T_{1-n}$ defines an injective morphism
$$F^\#\bar E_{1-n}/Ker (T.\bar\tau- \tau.F^\#T_{1-n})\hookrightarrow
E_{r-n},$$
we conclude that $T.\bar\tau- \tau.F^\#T_{1-n}=0$. We finish the proof
because for each $s\in Spec(R)$ we have $Im (T.\bar\tau- \tau.F^\#T_{1-n}))
\subset  m_s.  E_r $ for all $s\in
Spec(R)$ then, as $ E_r $ is locally free and $R$ is reduced,  we have
$T.\bar\tau- \tau.F^\#T_{1-n} =0$.

We repeat  the same argument to obtain $T_{l-n}:\bar E_{l-n}\to E_{r+l-n}$
with $0\leq l\leq n$.
\end{proof}

By looking at $\Upsilon$ and at the morphism defined in the last subsection,
$$\bar \beta: P_{p+D}^{n,ng+r }\to Cht_{p+D}^{n,ng+r}\times
Cht_{p+D}^{n,ng+r} \overset {\bar \Phi\times \bar
\Phi}\longrightarrow {\cal M}_{C}(n,ng+r,p+ D)\times {\cal M}_{C}(n,ng+r,p+
D),$$
we obtain, in the setting of elliptic sheaves,  the morphism $\bar \beta
.(\Upsilon\times \Upsilon)$, which    is the
analogous morphism to $\Phi\times \Phi$ defined in   subsection 3.3. We
again denote  this morphism   by $\bar \beta$.

From the last Lemma and Lemma \ref{isogenia} we obtain:

\begin{lem}\label{isogenia2}  If
$$[(\bar E_{j}, \bar i_{j},\bar \tau,g_{p+D}) ,(E_{j},i_{j},\tau,f_{p+D}) ]
\in \bar
\beta^{-1}(Z^{p+   D}_n)$$ defined over a field
$K$, then:

 1) There exists a $p+D$-isogeny between the elliptic sheaves with
$p+D$-level structures
$(\bar E_{j}, \bar i_{j},\bar \tau,g_{p+D})$ and $(E_{j},i_{j},\tau,f_{p+D})$.

or

2) There exists an isogeny between the elliptic sheaf  $(\bar E_j, \bar
i_j,\bar \tau )$  and the twisted
elliptic sheaf
 $(E_j,i_j,\tau ) \underset {\o_C} \otimes \o_C(-p-D)$.

or

3) $(F^{\#}\bar E_0)^\vee\underset {\o_C\otimes K}\otimes E_0 ((g+d)\infty
-p-D )\in \Theta_{n^2}$.
\end{lem}
\begin{proof}Let consider us the diagram  (not necessarily commutative)

$$ \xymatrix { \bar E_n   \ar[rr]^{T(1)} \ar[dr]      & &
   E_n(g+d-1)  \ar[dl]
\\ \bar E_1\ar[u]^{\bar l} &(\o_{p+ D})^n\otimes R &  E_1(g+d-1)\ar[u]^l
\\ F^{\#}\bar E_0\ar[u]^{\bar\tau }\ar[rr]^{F^{\#} T } \ar[ur] &  &F^{\#}
E_0(g+d-1)\ar[u]^{  \tau }
\ar[ul]    }
$$
$T$ is a morphism given by $[(\bar E_{j}, \bar i_{j},\bar \tau,g_{p+D})
,(E_{j},i_{j},\tau,f_{p+D}) ]
\in \bar
\beta^{-1}(Z^{p+   D}_n)$. The diagonal morphisms are the level structures,
and the vertical morphisms $\bar l$ and $l$
are induced by the "$i_j
$" and "$\bar i_j$". From Definition \ref{zeta2}, we have two possibilities:
either $T:\bar E_0  \to
 E_0 (g+d-1) $ takes values in  $E_0 (g+d-1)(-p-D)$ or the latter diagram
is commutative after
tensoring by $\o_{pD}$. In both cases, one obtains a morphism:
$$T(1).\bar l.\bar \tau.-l.\tau .F^{\#} T:F^{\#}\bar E_0  \to E_0 ((g+d
)\infty -p-D ).$$
 When $T.\bar l.\tau.t-\tau'. l'.t'.F^{\#} T=0$, we obtain
  cases 1) and 2) and when it is $\neq 0$ it tells us that $(F^{\#}\bar
E_0)^\vee\underset {\o_C\otimes K}\otimes
 E_0 ((g+d)\infty  -p-D)\in \Theta_{n^2}$.
\end{proof}
 Thanks to \cite{An2} Lemma 3.3, for $n=1$ it is possible to complete  3)
in terms of Drinfeld modules. It would be
very interesting to study this Lemma for rank-$n$-shtukas  and
rank-$n$-Drinfeld modules. For $C={\Bbb P}_1$,
because
$E_0$ and $\bar E_0$ are semi-stable  of $0$ degree
we have:
$$(F^{\#}\bar E_0)^\vee\underset {\o_C\otimes K}\otimes E_0 ((g+d)\infty
-p-D )\notin \Theta_{n^2}.$$
 One can
study
  explicit examples, as in    section 3, merely by changing monic
polynomials $"q(t)"$ for   polynomials $q(t)$ with
$q(\alpha_p)=1$, where
$m_p=(t-\alpha_p).{\Bbb F}_q[t]$.

\begin{rem} Given a  rank-$n$-elliptic sheaf $(   E_j,   i_j,  \tau )$, we
have that
$$h^0(E_{-n})=h^1(E_{-n})=0.$$
Thus,
$ E_{-n}$ is a semistable vector bundle. In this way,   there exists a
$d>>0$ such that for every pair of
elliptic sheaves $[(\bar E_{j}, \bar i_{j},\bar \tau ) ,(E_{j},i_{j},\tau
)]$ defined over a ring $R$ is
$$H^1(C\otimes R, {\bar  E_0}^\vee \underset{\o_C \otimes R } \otimes
E_0((g+d-1) )=0.$$
Therefore, $[(\bar E_{j}, \bar i_{j},\bar \tau ) ,(E_{j},i_{j},\tau )] \in
{\cal U}_d$.

In  the case of shtukas, this result is not true  unless   one considers
shtukas with truncation \cite{La1} and
\cite{La2}. In a following work, it would be interesting to  study   this
issue in the setting of
Lafforgue's shtuka compactifications \cite{La1}.

\end{rem}

 {\bf{ Acknowledgments}} I began to study this matter after a stay by
Prof.G.W.Anderson   in May 2002 at the University
of Salamanca. The conversations, suggestions and the encouragement of Prof.
G.W.Anderson have been invaluable to me for
studying this issue in the rank-$n$ case. I should like to express my
gratitude to Ricardo Alonso for his friendship and
company at all times and Jesus Muñoz Diaz and Laurent Lafforgue for their
friendship, desinterested help and Knidness
towards me.



\vskip1.5truecm { \'Alvarez V\'azquez, Arturo}\newline {\it
e-mail: } aalvarez@gugu.usal.es

\end{document}